\documentclass{article}

\usepackage[english]{babel}

\usepackage[a4paper,top=2cm,bottom=2cm,left=3cm,right=3cm,marginparwidth=1.75cm]{geometry}
\usepackage{amsmath}
\usepackage{amsthm}
\usepackage[colorlinks=true, allcolors=blue]{hyperref}
\usepackage{booktabs}
\usepackage{extarrows}
\usepackage{tikz-cd}
\usepackage{enumerate}
\usepackage{fancyhdr}
\usepackage{float}
\usepackage[nottoc]{tocbibind}
\usepackage{geometry}
\usepackage{indentfirst}
\usepackage{amssymb}
\usepackage{mathtools}
\usepackage{subfigure}
\usepackage{setspace}
\usepackage{tabularx}
\usepackage{graphicx}
\usepackage{booktabs}
\usepackage{multicol}
\usepackage{multirow}
\usepackage{longtable}
\usepackage{threeparttablex} 
\usepackage{booktabs}        
\usepackage{makecell}
\usepackage{tabularray}
\usepackage{algorithm}
\usepackage{comment}
\usepackage{listings}
\usepackage{xcolor}
\usepackage[subfigure]{tocloft}
\usepackage[format=hang,font=small,textfont=it]{caption}

\usepackage{titlesec}
\titlelabel{\thetitle.~}
\makeatletter
\renewenvironment{proof}[1][\proofname]{\par
  \pushQED{\qed}%
  \normalfont \topsep6\p@\@plus6\p@\relax
  \trivlist
  \item[\hskip\labelsep
        \bfseries
    #1\@addpunct{.}]\ignorespaces
}{%
  \popQED\endtrivlist\@endpefalse
}
\makeatother

\lstdefinestyle{macaulay}{
    basicstyle=\small\ttfamily,       
    backgroundcolor=\color{white},            
    frame=single,                             
    framerule=0.8pt,                          
    rulecolor=\color{black},                  
    tabsize=4,                                
    aboveskip=12pt,                           
    belowskip=12pt,                           
    breaklines=true,                          
    breakatwhitespace=true,                   
    showstringspaces=false,                   
    literate=                               
        {>}{{>\allowbreak}}{1}                
        {<}{{<}}{1}
        {:=}{{:=}}{2}
}

\DeclareCaptionFormat{listing}{\centering #1 #2 #3}
\newtheorem{thm}{Theorem}[section]
\newtheorem{prop}[thm]{Proposition}
\newtheorem{lemma}[thm]{Lemma}
\newtheorem{den}[thm]{Definition}
\newtheorem{coro}[thm]{Corollary}
\newtheorem{remark}[thm]{Remark}

\newtheorem{convention}[thm]{Convention}
\newtheorem{example}[thm]{Example}

\renewcommand\arraystretch{2}
\newcommand\mAut{\mathrm{Aut}}

\DeclareMathOperator{\mad}{ad}

\newcommand\mspan{\mathrm{Span}}
\newcommand\mnil{\mathrm{nil}}
\renewcommand\appendix{\par
    \setcounter{section}{0}
    \setcounter{subsection}{0}
    \gdef\thesection{\Alph{section}}}
\numberwithin{equation}{section}
\usepackage{hyperref}

\makeatletter
\renewcommand*\env@matrix[1][\arraystretch]{%
  \edef\arraystretch{#1}%
  \hskip -\arraycolsep
  \let\@ifnextchar\new@ifnextchar
  \array{*\c@MaxMatrixCols c}}
\makeatother

\makeatletter
\newcommand{\subjclass}[2][1991]{%
  \let\@oldtitle\@title%
  \gdef\@title{\@oldtitle\footnotetext{#1 \emph{Mathematics Subject Classification.} #2}}%
}
\newcommand{\keywords}[1]{%
  \let\@@oldtitle\@title%
  \gdef\@title{\@@oldtitle\footnotetext{\emph{Key words and phrases.} #1.}}%
}
\makeatother

\title{Classification of restricted Lie algebras of dimension 4}

\author{W. Liu, \; G.-S. Zhou }
\date{}

\newcommand{\Addresses}{{
  \bigskip
  \footnotesize

  W. Liu, \textsc{School of Mathematics and Statistics, Ningbo University, Ningbo 315211, China}\par\nopagebreak
  \textit{E-mail address}: \texttt{liuwei201029@gmail.com}

  \medskip

  G.-S. Zhou, \textsc{School of Mathematics and Statistics, Ningbo University, Ningbo 315211, China}\par\nopagebreak
  \textit{E-mail address}: \texttt{zhouguisong@nbu.edu.cn}
}}

	

\begin{document}

\maketitle

\begin{abstract}
Restricted Lie algebras of dimension up to $3$ over algebraically closed fields of positive characteristic  were classified by Wang and his collaborators in \cite{wang2013connected,wang2015classification}.
In this paper, we obtain a classification of restricted Lie algebras of dimension 4 over such fields.
\end{abstract}

\section{Introduction}

A \emph{restricted Lie algebra} (over a field of characteristic $p>0$) is a Lie algebra $L$ together with a self-map $[p]: L\to L$ satisfying certain conditions. Such a self-map $[p]$ is referred as a $p$-map or restricted map on $L$. In positive characteristic, Lie algebras which arise “in nature” are generally restricted. Classical examples include the associated Lie algebra of any associative algebra, the derivation algebra of any (Lie, associative, dtc.) algebra, the Lie algebra of any algebraic group, and the primitive space of any Hopf algebra. Restricted Lie algebras were first introduced and systematically studied by Jacobson \cite{jacobson1937derivation,jacobson1941restricted,jacobson1979lie}, and have played a predominant role in the theory of Lie algebras in positive characteristic, in connection with  classification of simple Lie algebras, algebraic groups, Hopf algebras and representation theory \cite{BLOCK1988115,milne,MR174052,strade2007lie,ModularLiealgebrasandresp}. A lot of studies have shown that in many respects restricted Lie algebras  bear a closer relation to Lie algebras of characteristic zero than ordinary Lie algebras of positive characteristic. 

The classification of Lie algebras of low dimension is a classical problem. Many classification results have appeared over various ground fields. Notably, there is a classification of Lie algebras of dimension up to $6$ over the complex number field and the real number field \cite{Mubarakzjanov1963classdim5,Mubarakzjanov1963solvdim6,Turkowski1990solvdim6,Snobl2014classificationsum}; of nilpotent Lie algebras of dimension up to $6$ over arbitrary ground fields \cite{Morozov1958nildim6,Gong1998thesis,Graaf2009nilcharno2,Graaf2012nil}; of solvable Lie algebras of dimension up to $4$ over arbitrary ground fields \cite{Graaf2004Classification}; and of non-solvable Lie algebras of dimension up to $6$ over finite fields as well as algebraically closed fields of positive characteristic \cite{strade2007lie}.  For higher dimensional cases, see \cite{classificationhistory,Le2023solvdim7} and the references therein. 

In comparison, less is known in literature about the classification of restricted Lie algebras of low dimension. These of dimension up to $3$ over algebraically closed fields were classified by Wang and his collaborators, as a part of their classification of connected Hopf algebras of low dimension in positive characteristic  \cite{wang2013connected,wang2015classification}. In addition, there is a classification of $p$-nilpotent restricted Lie algebras of dimension up to $4$ over perfect fields \cite{schneider2016classification}; and of $p$-nilpotent restricted Lie algebras of dimension  $5$ over algebraically closed fields of characteristic $p>3$ \cite{fivedimensionalpnil}. In \cite{5-dimensionalipnilincorrect}, there is given a list of $p$-nilpotent restricted Lie algebras of dimension $4$ over perfect fields of characteristic $p>3$. Unfortunately, this list is incomplete and redundant, as explained in \cite{fivedimensionalpnil}. 
Note that being $p$-nilpotent  is a very restrictive condition on restricted Lie algebras. It means that the $p$-map is nilpotent, which follows  necessarily that the Lie algebra itself is nilpotent.

In this paper, we achieve a classification of restricted Lie algebras of dimension $4$ over algebraically closed fields $\mathbb{F}$ of characteristic $p>0$ (see Theorem \ref{theclassification}). Naturally it yields a classification of primitively generated connected Hopf algebras  of dimension $p^4$ and  of  infinitesimal group schemes  of order $p^{4}$ and height $1$ over $\mathbb{F}$ (see  Remark \ref{Hopf-group-scheme}). By the classification results of \cite{wang2013connected,wang2015classification}, there are exactly  two isomorphism classes of restricted Lie algebras over $\mathbb{F}$ in dimension $1$; five isomorphism classes in dimension $2$; and $14$ (resp. $p+15$) isomorphism classes  for $p=2$ (resp. $p\geq 3$)  in dimension $3$.
Nevertheless, Theorem \ref{theclassification} tells us that  there are infinitely many isomorphism classes in dimension $4$. Actually, the isomorphism classes of our classification  consists of $2$ (resp. $3$) infinite one-parameter families, and a finite number of individuals or finite parameterized families for $p=2$ (resp. $p\geq 3$), as listed in Table \ref{mainresult}. The number of isomorphism classes that do not belong to the infinite one-parameter families are  explicitly described in Corollary \ref{numberisoclass}.

Now let us illustrate the method used in this paper. It is easy to see that two restricted Lie algebras $(L_1, [p]_1)$ and $(L_2, [p]_2)$ are isomorphic if and only if there is an isomorphism $g:L_{1}\rightarrow L_{2}$  of Lie algebras such that $g \circ [p]_{1}\circ g^{-1}$ and $[p]_{2}$ are conjugate  as $p$-maps on $L_2$. See  Definition \ref{conjugate} for the terminology of conjugate of two $p$-maps. So to classify restricted Lie algebras of dimension $d$ up to isomorphism, it suffices to finish the following two steps:
 \begin{itemize}
     \item Classify the Lie algebras of dimension $d$ up to isomorphism;
     \item Classify the $p$-maps up to conjugate on each one in the output of the first step.
 \end{itemize}
This is the approach that we take in this paper. It is similar to that of \cite{wang2013connected,wang2015classification,schneider2016classification}. As one can imagine, a lot of calculations have been taken in the course of this work. Though most of the calculations are straightforward, some of them are complicated. For those complicated (e.g. Remarks \ref{Groebner} and \ref{Groebner-2}), we have used the computer algebra system Magma \cite{magma}.

The organization of this paper is as follows. In Section \ref{pre}, we recall the basic definitions of restricted Lie algebras and review some basic facts. In Section \ref{Classify-Lie}, we provide a classification of Lie algebras of dimension $4$ over algebraically closed fields of positive characteristic. For these non-solvable, we simply follow the classification done by  Strade \cite{strade2007lie}. While for these solvable, our classification is an adjustment of the one given by de Graaf \cite{Graaf2004Classification}. This adjustment is necessarily because it is too complicated  to classify $p$-maps on some of the representative Lie algebras that occurred in the classification of de Graaf. Sections \ref{Classify-restrictedLie-1} - \ref{Classify-restrictedLie-3} are devoted to determine all restricted maps (up to conjugate) on each representative Lie algebra of dimension $4$ obtained in Section \ref{Classify-Lie}. The last section  presents our main result, offering the classification of restricted Lie algebras of dimension $4$ over algebraically closed fields of positive characteristic. 

\paragraph{Notation and Conventions.} Throughout the paper, we work over a fixed  field $\mathbb{F}$. It is assumed to be algebraically closed of positive characteristic $p$, unless otherwise stated. All algebras and vector spaces are taken over $\mathbb{F}$.  The set of nonzero elements of $\mathbb{F}$  is denoted by $\mathbb{F}^{\times}$. For any $\alpha \in \mathbb{F}^{\times}$, the equation $x^{p-1}-\alpha=0$ has $p-1$ distinct roots in $\mathbb{F}$. We fix one and denote it by $\alpha^{\frac{1}{p-1}}$.  
 
The prime field of characteristic $p$ is denoted by $\mathbb{F}_{p}$.  The set of nonzero elements of $\mathbb{F}_{p}$ is denoted by $\mathbb{F}_{p}^{\times}$.  We need specific subsets of $\mathbb{F}_{p}^{\times}$ to parameterize certain families of the classifying restricted Lie algebras. Since $\mathbb{F}_{p}^{\times}$ is a cyclic group, we may choose a generator $\xi_{p} \in \mathbb{F}_{p}^{\times}$ and set $\Xi_{p}:=\{\xi_{p}^{r} ~ |~ 0\leq r\leq \frac{p-1}{2}\}$. Note that $\Xi_{p}\cup \Xi_{p}^{-1} = \mathbb{F}_{p}^{\times}$ and $\Xi_{p}\cap \Xi_{p}^{-1} = \{\pm1\}$. In addition, let $Q_{p}:=\{ \alpha\in \mathbb{F}_{p}^{\times} ~|~ \alpha^{\frac{p-1}{2}}=1 \}$. It is the set of quadratic residues modulo $p$.

The center of a Lie algebra $L$ is denoted by $Z(L)$. For a Lie algebra $L$ of dimension $4$ which is specified with a  basis $x,y,z,w$, we denote by  $\sigma_{a,b,c,d}:L\to L$  for any $a,b,c,d \in \mathbb{F}$ the linear transformation given  by $x\mapsto ax,  y\mapsto by,  z\mapsto cz,  w\mapsto dw $.

\section{Preliminaries}\label{pre}

In this section, we recall basic notions of restricted Lie algebras that should be used in the sequel. We refer to \cite{jacobson1979lie, milne, ModularLiealgebrasandresp} for basic references. Throughout this section, the base field $\mathbb{F}$ is only assumed to be of characteristic $p>0$. It is not necessarily algebraically closed.

Let $L$ be a Lie algebra over $\mathbb{F}$. We denote by $\mad=\mad_L: L\to \mathfrak{gl}(L)$ the adjoint representation. So $\mad(x)(y) =[x,y]$ for any $x, y\in L$. For $0<i<p$ and $x_0,x_1\in L$, we set
\[
s_i(x_0,x_1):= -\frac{1}{i} \sum_{u} \mad(x_{u(1)})\mad(x_{u(2)})\cdots \mad(x_{u(p-1)}) (x_{1}),
\]
where $u$ runs over the maps $\{1, 2, \ldots, p-1\} \to \{0,1\}$ taking $i$ times the value $0$. 
Let $L[T]$ be the space of polynomials in the variable $T$ with coefficients in $L$. Consider it as a Lie algebra over the polynomial ring $\mathbb{F}[T]$ by base change. Then  one has 
\begin{eqnarray}\label{adT-1}
\big(\mad_{L[T]}(x_0T +x_1)\big)^{p-1}(x_0) =  \sum_{i=1}^{p-1} is_{i}(x_0,x_1) T^{i-1}.
\end{eqnarray}
It is a useful formula to calculate $s_i(x_0,x_1)$ in practice.

\begin{den}\label{pmapdef}
A map $[p]:L\to L, ~ x\mapsto x^{[p]}$, is said to be  \emph{$p$-semilinear} if it satisfies 
 \begin{enumerate}
    \item $(\alpha x)^{[p]}=\alpha^{p}x^{[p]}$ for any $x\in L$ and $\alpha\in \mathbb{F}$;
    \item $(x+y)^{[p]}=x^{[p]}+y^{[p]}+\sum\limits_{i=1}^{p-1}s_{i}(x,y)$ for any $x, y\in L$.
\end{enumerate}
A \emph{$p$-map} on $L$ is a $p$-semilinear map such that $\mad x^{[p]} = (\mad x)^{p}$ for any $x\in L$.
\end{den}

\begin{example}
Given an associative algebra $A$, one may consider it as a Lie algebra under the natural bracket $[x,y]:=xy-yx$. It is well-known that the $p$-th power map  $x\mapsto x^p$ is a $p$-map on $A$. See \cite[Proposition 10.38]{milne} for a reference.     
\end{example}

The next result provides a method to describe or construct $p$-maps on an arbitrary Lie algebra. It should be well-known to experts, but we failed to find a reference. See \cite[Theorem 2.3]{ModularLiealgebrasandresp} for a related but difference result. We provide a detailed proof for the reader's convenience. 

\begin{lemma}\label{pmapconstruct}
Let $L$ be a Lie algebra and $\{x_j\}_{j\in J}$ a basis of $L$. Then every family $\{f_j\}_{j\in J}$ of elements in $L$ uniquely determines a $p$-semilinear map $[p]$ on $L$ such that $x_j^{[p]} =f_j$ for any $j\in J$. Moreover, the map $[p]$ is a $p$-map if and only if $\mad f_j = (\mad x_j)^p$ for each $j\in J$.
\end{lemma}

\begin{proof}
The uniqueness of the first statement is clear. To see the existence, we consider $L$ as a subspace of the universal enveloping algebra $U(L)$ in the natural way. For every family  $\{f_j\}_{j\in J}$ of elements in $L$, one may define a self-map $[p]$ on $L$ by
 \[
 x \mapsto x^{[p]}:= x^p - \sum_{j\in J} a_j^{p} x_j^p + \sum_{j\in J} a_j^{p} f_j,
 \]
where $a_j\in \mathbb{F}$ is determined by $x=\sum_{j\in J} a_j x_j$. To see that it is indeed a self-map on $L$, it suffices to show $x^p - \sum_{j\in J} a_j^{p} x_j^p\in L$. To this end, let us expresses $x$ as a finite sum $x = \sum_{r=1}^n a_{r}x_{j_r}$. Since the $p$-th power map on $U(L)$ is a $p$-map, it follows that 
\begin{eqnarray*}
x^p - \sum_{r=1}^n a_{r}^px_{j_r}^p  = \sum_{i=1}^{p-1}  \sum_{r=2}^{n} s_i\left(\sum_{s=1}^{r-1} a_{s} x_{j_s},  a_{r} x_{j_r}\right) \quad  \in \quad  L. 
\end{eqnarray*}
Next we show that the map $[p]$ is a $p$-semilinear map on $L$. Clearly, $(ax)^{[p]} = a^px^{[p]}$ for any $a\in \mathbb{F}$ and $x\in L$. In addition, for any  $x, y\in L$ one has
\[
(x+y)^{[p]} - x^{[p]} -y^{[p]} = (x+y)^p - x^p -y^p = \sum_{i=1}^{p-1} s_i(x, y).
\]
Thus the map $[p]$ is indeed a $p$-semilinear map satisfies that $x_j^{[p]} =f_j$ for any $j\in J$.

To see the second statement, it suffices to show that for any $p$-semilinear map $[p]$ on $L$, the subset $L_0:= \{x\in L~ |~ \mad x^{[p]} = (\mad x)^p\}$ is a subspace of $L$. Readily, $ax\in L_0$ for any $a\in \mathbb{F}$ and $x\in L_0$. In addition, if $x, y\in L_0$ then $x+y \in L_0$. It is because
\begin{eqnarray*}
\mad (x+y)^{[p]} &=& \mad \left(x^{[p]} + y^{[p]} +\sum_{i=1}^{r-1} s_i(x,y)\right) \\
               &=& (\mad x)^{p} + (\mad y)^{p} + \sum_{i=1}^{r-1} s_i (\mad x, \mad y) \\
               &=& (\mad x+ \mad y)^p = (\mad (x+y))^p.
\end{eqnarray*}
Here, we used the fact that $\mad: L\to \mathfrak{gl}(L)$ is a homomorphism of Lie algebras and the $p$-th power map on $\mathfrak{gl}(L)$ is a $p$-map. Thus, $L_0$ is a subspace of $L$, as desired.
\end{proof}

\begin{convention}
Due to Lemma \ref{pmapconstruct},  $p$-semilinear maps (particularly, $p$-maps)  on a Lie algebra $L$ will be described in the sequel  by a list of the form  
$$x_j^{[p]} =f_j, \quad j\in J,$$ where $\{x_j\}_{j\in J}$ is a basis of $L$ and $f_j\in L$. Frequently, we only write nonzero ones.
\end{convention}

\begin{den}
A \emph{restricted Lie algebra} is a Lie algebra together with a $p$-map on it. Precisely, it is a pair $(L, [p])$, where $L$ is a Lie algebra and $[p]$ is a $p$-map on $L$.
\end{den}

\begin{den}
Let $(L_1, [p]_1)$ and $(L_2, [p]_2)$ be restricted Lie algebras. A \emph{homomorphism} from $(L_1, [p]_1)$ to $(L_2, [p]_2)$ is a homomorphism $g:L_1 \to L_2$ of Lie algebras such that $$g\circ [p]_1= [p]_2\circ g.$$ 
If in addition $g$ is invertible, then it is called an \emph{isomorphism}. We say that $(L_1, [p]_1)$ and $(L_2, [p]_2)$ are \emph{isomorphic} if there is an isomorphism from $(L_1, [p]_1)$ to $(L_2, [p]_2)$.    
\end{den}

Let $\mAut(L)$ be the group of automorphisms of $L$ as Lie algebras. Then for any $\varphi\in \mAut(L)$, if $[p]$ is a $p$-semilinear map (resp. $p$-map) on $L$ then so is $\varphi\circ [p] \circ \varphi^{-1}$. In this way, one obtains a group action of $\mAut(L)$ on the set of a $p$-semilinear maps (resp. $p$-maps) on $L$.

\begin{den}\label{conjugate}
We say that two $p$-semilinear maps $[p]_1, [p]_2$ on a Lie algebra $L$ are \emph{conjugate} (or \emph{equivalent}) if there exists $\varphi\in \mAut (L)$ such that $\varphi\circ [p]_1\circ \varphi^{-1} =[p]_2$.
\end{den}

To conclude this section,  we give several conjugate invariances for $p$-maps on a Lie algebra $L$.  For any $p$-map $[p]$ on $L$, any subspace $V$ of $L$ and any integer $r\geq 0$, we set  
\[
V^{[p]^{r}}:= \mspan\{v^{[p]^{r}} ~|~ v\in V \}.
\]

\begin{lemma}\label{ConjInvar}
Let $[p]_1$ and $[p]_2$ be conjugate $p$-maps on a Lie algebra $L$. Then $\dim L^{[p]_1^r} = \dim L^{[p]_2^r}$, $\dim Z(L)^{[p]_1^r} = \dim Z(L)^{[p]_2^r}$ and $\dim [L,L]^{[p]_1^r} = \dim [L,L]^{[p]_2^r}$ for any $r\geq1$.
\end{lemma}
\begin{proof}
It is easy to see.
\end{proof}


\section{Classification of Lie algebras of dimension $4$}
\label{Classify-Lie}

This section is devoted to present a complete representatives of Lie algebras of dimension $4$. The classification problem is split into the following three cases: (1) solvable with a $3$-dimensional abelian ideal; (2) solvable without $3$-dimensional abelian ideals; and (3) non-solvable. For the last case, we simply follow the classification that done by Strade \cite{strade2007lie}; while for the first two cases, our classification is an adjustment of the one given by de Graaf \cite{Graaf2004Classification}.

The following Lemma is employed to decide whether two of classifying Lie algebras are isomorphic. Recall that every finite dimensional Lie algebra $L$ has a unique maximal nilpotent ideal. It  is called the nilradical of $L$, and is denoted by $\mnil(L)$.

\begin{lemma}\label{iso-criterion}
Let $L$ and $M$ be solvable Lie algebras of dimension $n$. Assume that $\mnil(L)$ and $\mnil(M)$ are both abelian and of dimension $n-1$. Let $x\in L\backslash \mnil(L)$ and $y\in M\backslash \mnil(M)$. Let $A$ be the matrix of $\mad_{L}(x)|_{\mnil(L)}$ with respect to any chosen basis of  $\mnil(L)$.  Let $B$ be the matrix of $\mad_{M}(y)|_{\mnil(M)}$ with respect to any chosen basis of  $\mnil(M)$. Then $L$ and $M$ are isomorphic as Lie algebras if and only if $A$ and $kB$ are conjugate as matrices for some $k\in \mathbb{F}^{\times}$.
\end{lemma}

\begin{proof}
Let $\varphi:L\to M$ be an isomorphism of Lie algebras. Then $\varphi(x)\in ky+\mnil(M)$ for some $k\in \mathbb{F}^{\times}$, and $\varphi$ restricts to a vector space isomorphism  $\varphi|_{\mnil(L)}:\mnil(L)\to \mnil(M)$. Thus, $\varphi|_{\mnil(L)}$ specifies an invertible matrix $P\in \mathrm{GL}_{n-1}(\mathbb{F})$ with respect to the given bases for $\mnil(L)$ and $\mnil(M)$. Since $\varphi|_{\mnil(L)}\circ\mad_{L}(x)|_{\mnil(L)}= \mad_{M}(ky)|_{\mnil(M)}\circ\varphi|_{\mnil(L)}$, one obtains $PA P^{-1} = kB$.

Conversely, assume $PAP^{-1}=kB$ for some $P\in \mathrm{GL}_{n-1}(\mathbb{F})$ and $k\in\mathbb{F}^{\times}$. Then there exists a unique linear  isomorphism $\psi:\mnil(L)\to \mnil(M)$ whose matrix is $P$, with respect to the given bases of $\mnil(L)$ and $\mnil(M)$. 
Let  $\tilde{\psi}:L\to M$ be the linear map given by $\tilde{\psi}(x)=ky$ and $\tilde{\psi}|_{\mnil{(L)}}=\psi$. It is easy to check that $\tilde{\psi}$ is an isomorphism of Lie algebras.
\end{proof}

In order to parameterize the classifying Lie algebras in a compact form, we utilize the group action of the permutation group $S_{3}$ on $\mathbb{F}^{\times}\times\mathbb{F}^{\times}$ given by
\begin{eqnarray}\label{actionS3}
  (12)\cdot(\xi,\eta)=(\xi^{-1},\xi^{-1}\eta), \quad \quad  (23)\cdot(\xi,\eta)=(\eta,\xi).
\end{eqnarray}

\begin{prop}\label{Liewith3ideal}
Every solvable Lie algebra of dimension $4$ with a $3$-dimensional abelian ideal is isomorphic to one of the following: 
    \begin{enumerate}
        \item $L_{1}=abelian$;
        \item $L_{2}=\langle ~x,y,z,w~|~[w,x]=y~\rangle$; 
        \item $L_{3}=\langle ~x,y,z,w~|~[w,x]=y,~[w,y]=z~\rangle$; 
        \item $L_{4}=\langle ~x,y,z,w~|~[w,x]=x~\rangle$; 
        \item $L_{5}(\xi)=\langle ~x,y,z,w~|~[w,x]=x,~[w,y]=\xi y~\rangle, $ where $\xi\in \mathbb{F}^{\times}$;
        \item $L_{6}(\xi,\eta)=\langle ~x,y,z,w~|~[w,x]=x,~[w,y]=\xi y,~[w,z]=\eta z~\rangle, $ where $\xi,\eta\in \mathbb{F}^{\times}$;
        \item $L_{7}=\langle ~x,y,z,w~|~[w,x]=y,~[w,z]=z~\rangle$;
        \item $L_{8}(\xi)=\langle ~x,y,z,w~|~[w,x]=x+y,~[w,y]=y,~[w,z]=\xi z~\rangle,$ where $\xi\in \mathbb{F}$;
        \item $L_{9}=\langle ~x,y,z,w~|~[w,x]=x+y,~[w,y]=y,~[w,z]=x+z~\rangle$.
    \end{enumerate}
    Moreover, these Lie algebras are mutually non-isomorphic except the following cases: (1)  $L_{5}(\xi)\cong L_{5}(\xi^{-1})$ for any $\xi\in\mathbb{F}^{\times}$; (2) $L_{6}(\xi,\eta)\cong L_{6}(\tau\cdot (\xi, \eta))$  for any $\tau \in S_{3}$ and $\xi,\eta \in \mathbb{F}^{\times}$.
\end{prop}
\begin{proof}
Let $L$ be a solvable Lie algebra of dimension $4$ with a $3$-dimensional abelian ideal, say $I$. Fix an element  $w\in L\backslash I$. It is clear that the Lie structure of $L$ is determined by $\mad w: I\to I$. Since $\mathbb{F}$ is algebraically closed, we may choose a basis $x, y, z$ of $I$ such that the matrix of $\mad w$ with respect to it  is of one of the following  forms
    \begin{equation*}
        \begin{pmatrix}[1.3]
            \lambda_{1} & 0 & 0\\
            0 & \lambda_{2} & 0\\
            0 & 0 & \lambda_{3}
        \end{pmatrix}
        ,\quad
        \begin{pmatrix}[1.3]
            \lambda_{1} & 0 & 0\\
            1 & \lambda_{1} & 0\\
            0 & 0 & \lambda_{2}
        \end{pmatrix}
        ,\quad
        \begin{pmatrix}[1.3]
            \lambda_{1} & 0 & 0\\
            1 & \lambda_{1} & 0\\
            0 & 1 & \lambda_{1}
        \end{pmatrix}
    \end{equation*}
   where $\lambda_{1},\lambda_{2},\lambda_{3}\in\mathbb{F}$.   
   Firstly consider that $\mad w$ is of the left-most form with respect to $x, y,z$.  In this case, we may assume that if $\lambda_{i}=0$ then $\lambda_{i+1}=0$ (by changing the order of $x, y,z$).
    \begin{enumerate}
        \item If $\lambda_{1}=0$ then we obtain a presentation of $L$ as $L_{1}$.
        \item If $\lambda_{1}\neq0$ and $\lambda_{2}=0$ then we obtain a presentation of $L$ as $L_{4}$  by replacing $w$ with $\lambda_{1}^{-1}w$.
        \item If $\lambda_{1},\lambda_{2}\neq0$ and $\lambda_{3}=0$ (resp. $\lambda_{3}\neq 0$) then we obtain a presentation of $L$ as $L_{5}(\lambda_{2}/\lambda_{1})$  (resp. $L_{6}(\lambda_{2}/\lambda_{1},\lambda_{3}/\lambda_{1})$)  by replacing $w$ with $\lambda_{1}^{-1}w$.
    \end{enumerate}
    Next assume $\mad w$ is of the middle form with respect to $x, y,z$.
    \begin{enumerate}
        \item If $\lambda_{1}=\lambda_{2}=0$ then we obtain a presentation of $L$ as $L_{2}$.
        \item If $\lambda_{1}=0$ and $\lambda_{2}\neq0$ then we obtain a presentation of $L$ as $L_{7}$ by replacing $w$ with $\lambda_{2}^{-1}w$;
        \item If $\lambda_{1}\neq0$  then we obtain a presentation of $L$ as $L_{8}(\lambda_{2}/\lambda_{1})$  by replacing $y$ and $w$ with $\lambda_{1}^{-1}y$ $\lambda_{1}^{-1}w$ respectively.
    \end{enumerate}
    Finally  assume $\mad w$ is of the right-most form with respect to $x, y,z$.
    \begin{enumerate}
        \item If $\lambda_{1}=0$ then  we obtain a presentation of $L$ as $L_{3}$.
        \item If $\lambda_{1}\neq0$ then  we obtain a presentation of $L$ as $L_{9}$ by replacing $y$, $z$, $w$ with $\lambda_{1}^{-2}z$, $\lambda_{1}^{-1}y$ $\lambda_{1}^{-1}w$ respectively.
    \end{enumerate}
In conclusion, $L$ is isomorphic to one of the Lie algebras that listed in the Proposition.

Now we show the last statement. The listed Lie algebras are considered to share the same underling vector space with a basis $x,y,z,w$. Let  $L_{5}=L_{5}(\xi), L_{6} = L_{6}(\xi,\eta)$ and $L_{8} =L_{8}(\xi)$ for simplicity. Note that $L_{1}$, $L_{2}$, $L_{3}$ are nilpotent, while the others are not. In fact, the nilradical $\mnil(L_{i})$ is abelian and spanned by $x, y, z$ for $i=4,\ldots, 9$. Therefore, $L_{i} \not\cong L_{j}$ for $1\leq i\leq 3$ and $4\leq j\leq 9$. Moreover, $L_{1}$, $L_{2}$, $L_{3}$ are mutually non-isomorphic because their center are of different dimension. According to the matrix of $\mad_{L_{i}}(w)|_{\mnil(L_{i})}$ with respect to $x,y,z$, Lemma \ref{iso-criterion} tells us readily that $L_{i}\neq L_{j}$ for $4\leq i<j\leq 9$,  $L_{8}(\xi) \neq L_{8} (\xi')$ for $\xi\neq \xi'$, $L_{5}(\xi) \cong L_{5} (\xi')$ if and only if  $\xi^{-1}= \xi'$ for $\xi\neq \xi'$, and $L_{6}(\xi,\eta)\cong L_{6}(\xi',\eta')$ if and only if $(\xi,\eta)=\tau \cdot (\xi',\eta')$ for some $\tau \in S_{3}$.
\end{proof}

\begin{remark}
Solvable Lie algebras of dimension $4$ over a field of arbitrary characteristic have been classified in \cite[Section 5]{Graaf2004Classification}. The representatives presented there are labeled by $$M^{1},~ M^{2},~ M^{3}_{a},~ M^{4},~ M^{5}, ~M^{6}_{a,b},  ~ M^{7}_{a,b}, ~ M^{8}, ~ M^{9}_{a}, ~ M^{10}_{a}, ~ M^{11}_{a,b}, ~ M^{12}, ~ M^{13}_{a},~ M^{14}_{a},$$ where $a, b$ are scalars that satisfy certain conditions. It is easy to see that the Lie algebras in the first seven families have a $3$-dimensional abelian ideal, while the others do not. However, the representatives $M^{1}, M^{2}, M^{3}_{a}, M^{4}, M^{5}, M^{6}_{a,b}, M^{7}_{a,b}$ are different from those given in Proposition \ref{Liewith3ideal}. The difference is due to that we use the Jordan normal form to describe $\mad w$, while in  \cite{Graaf2004Classification}  the author utilized the rational normal form. It turns out that it is much easier to classify the $p$-maps up to conjugate on the representative Lie algebras $L_{i}$ than that on $M^{j}$.
\end{remark}

\begin{prop}\label{Liewithout3ideal}
  The isomorphism classes of solvable Lie algebras of dimension 4 without 3-dimensional abelian ideals are as follows:
    \begin{enumerate}
        \item $N_{1}=\langle ~x,y,z,w~|~[y,x]=x,~[w,z]=z~\rangle; $
        \item $N_{2}=\langle ~x,y,z,w~|~[z,x]=y,~[w,x]=x,~[w,y]=2y,~[w,z]=z~\rangle; $
        \item $N_{3}(\xi)=\langle ~x,y,z,w~|~[z,x]=y,~[w,x]=x+\xi z,~[w,y]=y,~[w,z]=x~\rangle, $ where $\xi\in\mathbb{F}^{\times};$
        \item $N_{4}=\langle ~x,y,z,w~|~[z,x]=y,~[w,x]=z,~[w,z]=x~\rangle; $
        \item $N_{5}=\langle ~x,y,z,w~|~[z,x]=y,~[z,y]=x,~[w,x]=x,~[w,z]=z~\rangle, $ where $p=2$.
    \end{enumerate}
\end{prop}
\begin{proof}
By \cite[Section 5]{Graaf2004Classification}, we already have a complete representative of solvable Lie algebras of dimension $4$ without 3-dimensional abelian ideals over a field of arbitrary characteristic. They are labeled by $M^{8},  M^{9}_{a}, M^{10}_{a}, M^{11}_{a,b},  M^{12},  M^{13}_{a}, M^{14}_{a},$ where $a, b$ are scalars that satisfy certain conditions. If the base field is assumed to be algebraically closed, then one may readily see that the parameter sets of
$M^{9}_{a}$ and $M^{10}_{a}$ are empty and there exist the following correspondences:
    \[
    M^{8}\cong N_{1},~ M_{a,b}^{11}\cong M_{1,0}^{11}\cong N_{5},~ M^{12}\cong N_{2},~ M^{13}_{a}\cong N_{3}(\xi),~ M^{14}_{a}\cong M^{14}_{1}\cong N_{4}.
    \]
The result follows directly.
\end{proof}

\begin{remark}\label{N5p3}
The defining relations of $N_{5}$ fail to yield a Lie algebra for $p\geq 3$ because 
    \[
    [w,[x,z]]+[z,[w,x]]+[x,[z,w]]=2y\neq 0.
    \]
\end{remark}

\begin{prop}\label{nonsolvLie}
The isomorphism classes of non-solvable Lie algebras of dimension 4 are as follows:
    \begin{enumerate}
        \item $\mathfrak{gl}_{2}=\langle ~x,y,z,w~|~[y,x]=-z,~[z,x]=2x,~[z,y]=-2y~\rangle$, where $p\geq 3$;
        \item $W_{1}=\langle ~x,y,z,w~|~[y,x]=x,~[z,x]=y,~[z,y]=z~\rangle, $ where $p=2$;
        \item $W_{2}=\langle ~x,y,z,w~|~[y,x]=x,~[z,x]=y,~[z,y]=z,~[w,x]=z~\rangle, $ where $p=2$.
    \end{enumerate}
\end{prop}
\begin{proof}
It is {\cite[Theorem 4.1]{strade2007lie}}. Though the base field is assumed to be finite in \cite{strade2007lie}, a similar argument applies in our context that the base field is algebraically closed.
\end{proof}

\begin{remark}
The defining relations of $\mathfrak{gl}_{2}$  give rise to a Lie algebra for arbitrary $p$. Furthermore, $\mathfrak{gl}_{2}$ is exactly the general lineal Lie algebra $\mathfrak{gl}_{2}(\mathbb{F})$ via the identification    
\[
 x=\begin{pmatrix}[1.3]
     0 &1~ \\
     0 & 0~\\
 \end{pmatrix},~
 y=\begin{pmatrix}[1.3]
     0 & 0~\\
     1& 0~\\
 \end{pmatrix},~
 z=
 \begin{pmatrix}[1.3]
     1 & 0\\
     0 & -1~
 \end{pmatrix},~
 w=\begin{pmatrix}[1.3]
     1 & 0~\\
     0 & 1~
 \end{pmatrix}.
 \]
Though $\mathfrak{gl}_{2}$ is not solvable for $p\geq 3$, it is for $p=2$. Actually, there is an isomorphism $\mathfrak{gl}_{2} \xrightarrow{\cong} L_{2}$ of Lie algebras given by $x\mapsto w, y\mapsto x , z\mapsto y, w\mapsto z$ when $p=2$.  
\end{remark}

\begin{remark}
 The defining relations of $W_{1}$ also give rise to a Lie algebra for $p\geq 3$. In fact, there is an isomorphism  $W_{1} \xrightarrow{\cong} \mathfrak{gl}_{2}$ of Lie algebras given by $x\mapsto x, y\mapsto \frac{1}{2}z , z\mapsto -\frac{1}{2}y, w\mapsto w$. However, the defining relations of $W_{2}$ fail to yield a Lie algebra for $p\geq 3$ because 
    \[
    [w,[y,x]]+[y,[x,w]]+[x,[w,y]]=2z\neq 0.
    \]
\end{remark}

\section{The $p$-maps on solvable Lie algebras of dimension $4$ with a $3$-dimensional abelian ideal}
\label{Classify-restrictedLie-1}

In this section, we classify up to conjugate the $p$-maps  on the Lie algebras that presented in Proposition \ref{Liewith3ideal}. Note that these Lie algebras form a full list of representatives  of solvable Lie algebras of dimension $4$ with a $3$-dimensional abelian ideal. 

\subsection{Case $L_{1}$}

Firstly, we classify the $p$-maps on $L_{1}$, which is the abelian Lie algebra of dimension $4$. 

\begin{prop}
   The conjugacy classes of the $p$-maps on $L_{1}$ are as follows:
    \begin{enumerate}
        \item Trivial $p$-map;
        \item $x^{[p]}=y,~y^{[p]}=0,~z^{[p]}=0,~w^{[p]}=0$;
        \item $x^{[p]}=y,~y^{[p]}=0,~z^{[p]}=w,~w^{[p]}=0$;
        \item $x^{[p]}=y,~y^{[p]}=z,~z^{[p]}=0,~w^{[p]}=0$;
        \item $x^{[p]}=y,~y^{[p]}=z,~z^{[p]}=w,~w^{[p]}=0$;
        \item $x^{[p]}=x,~y^{[p]}=0,~z^{[p]}=0,~w^{[p]}=0$;
        \item $x^{[p]}=x,~y^{[p]}=z,~z^{[p]}=0,~w^{[p]}=0$;
        \item $x^{[p]}=x,~y^{[p]}=y,~z^{[p]}=0,~w^{[p]}=0$;
        \item $x^{[p]}=x,~y^{[p]}=z,~z^{[p]}=w,~w^{[p]}=0$;
        \item $x^{[p]}=x,~y^{[p]}=y,~z^{[p]}=w,~w^{[p]}=0$;
        \item $x^{[p]}=x,~y^{[p]}=y,~z^{[p]}=z,~w^{[p]}=0$;
        \item $x^{[p]}=x,~y^{[p]}=y,~z^{[p]}=z,~w^{[p]}=w$.
    \end{enumerate}    
The $p$-maps in (1) to (5) are nilpotent, while the others are not.
\end{prop}
\begin{proof}
    By \cite[Proposition A.1]{wang2015classification}, the list follows. Moreover, they don't conjugate to each other since they have different 3-tuple of invarients: $\dim L_{1}^{[p]}$, $\dim L_{1}^{[p]^{2}}$, $\dim L_{1}^{[p]^{3}}$.
\end{proof}

\subsection{Case $L_{2}$}

Recall that the Lie algebra $L_{2}$ has a basis $x,y,z,w$ with one structure relation $[w,x]=y.$ Clearly, $Z(L_{2})=\mspan\{y,z\}$. In addition, a linear map $\varphi: L_{2} \to L_{2}$ is a homomorphism of Lie algebras if and only if  $\varphi(y), \varphi(z) \in \mspan\{y,z\}$ and $[\varphi(w), \varphi(x)] = \varphi(y)$. It follows that $\varphi \in \mAut(L_{2})$ if and only if it is given by 
\begin{equation}\label{AutoLb}
    \begin{aligned}
    x&\mapsto a_{1}x+a_{2}y+a_{3}z+a_{4}w, \quad  \quad \quad&
    y&\mapsto (a_{1}d_{4}-a_{4}d_{1})y,\\
    z&\mapsto c_{2}y+c_{3}z, \quad \quad \quad & 
    w&\mapsto d_{1}x+d_{2}y+d_{3}z+d_{4}w
    \end{aligned}
    \end{equation}
for some $a_i, c_i, d_i\in \mathbb{F}$, which satisfy $(a_{1}d_{4}-a_{4}d_{1}) c_{3}\neq 0$.

It is easy to verify that $(\mad x)^{p} = (\mad y)^{p} = (\mad z)^{p} = (\mad w)^{p} =0$.  Thus by Lemma \ref{pmapconstruct}, the $p$-maps on $L_{2}$ are the $p$-semilinear maps of the form:
\begin{eqnarray}
x^{[p]}=\alpha_{1}y+\alpha_{2}z, \quad y^{[p]}=\beta_{1}y+\beta_{2}z, \quad z^{[p]}=\gamma_{1}y+\gamma_{2}z, \quad w^{[p]}=\delta_{1}y+\delta_{2}z,
\end{eqnarray}
where $\alpha_i, \beta_i, \gamma_i, \delta_i\in \mathbb{F}$. The result below classifies the $p$-maps on $L_{2}$.

\begin{prop}\label{Lbp}
The conjugacy classes of the $p$-maps on $L_{2}$ are as follows:
    \begin{enumerate}
        \item Trivial $p$-map, for $p\geq 3$;
        \item $x^{[p]}=0,~y^{[p]}=0,~z^{[p]}=0,~w^{[p]}=y$;
        \item $x^{[p]}=0,~y^{[p]}=0,~z^{[p]}=0,~w^{[p]}=z$;
        \item $x^{[p]}=z,~y^{[p]}=0,~z^{[p]}=0,~w^{[p]}=y$;
        \item $x^{[p]}=0,~y^{[p]}=z,~z^{[p]}=0,~w^{[p]}=0$, for $p\geq 3$;
        \item $x^{[p]}=y,~y^{[p]}=z,~z^{[p]}=0,~w^{[p]}=0$;
        \item $x^{[p]}=0,~y^{[p]}=0,~z^{[p]}=y,~w^{[p]}=0$;
        \item $x^{[p]}=z,~y^{[p]}=0,~z^{[p]}=y,~w^{[p]}=0$;
        \item $x^{[p]}=0,~y^{[p]}=y,~z^{[p]}=0,~w^{[p]}=0$;
        \item $x^{[p]}=z,~y^{[p]}=y,~z^{[p]}=0,~w^{[p]}=0$;
        \item $x^{[p]}=0,~y^{[p]}=y+z,~z^{[p]}=0,~w^{[p]}=0$, for $p\geq 3$;
        \item $x^{[p]}=z,~y^{[p]}=y+z,~z^{[p]}=0,~w^{[p]}=0$;
        \item $x^{[p]}=0,~y^{[p]}=0,~z^{[p]}=z,~w^{[p]}=0$, for $p\geq 3$;
        \item $x^{[p]}=y,~y^{[p]}=0,~z^{[p]}=z,~w^{[p]}=0$;
        \item $x^{[p]}=0,~y^{[p]}=y+\lambda z,~ z^{[p]}=z,~w^{[p]}=0$,  $\lambda\in \mathbb{F}$. Here, $\lambda_{1}, \lambda_{2}\in \mathbb{F}$ represent conjugate $p$-maps if and only if $\lambda_{1}^{p(p-1)}(\lambda_{2}^{p-1}+1)^{p+1} =\lambda_{2}^{p(p-1)}(\lambda_{1}^{p-1}+1)^{p+1}.$
    \end{enumerate}
    The $p$-maps in (1) to (8) are nilpotent, while the others are not.
\end{prop}

The proof of this proposition will be addressed at the end of this subsection.

\begin{remark}
Let $[p]$ be the $p$-map on $L_{2}$ given by $x^{[p]} =0$, $y^{[p]} =0$, $z^{[p]} =0 $ and $w^{[p]}=0$. It is the trivial $p$-map when $p\geq 3$. However, it is not for $p=2$. Indeed, one has 
\[
(x+w)^{[2]} =x^{[2]} +w^{[2]} + [w,x] = y.
\]
Actually, $\varphi\circ [2]\circ \varphi^{-1}$ is the one given in Proposition \ref{Lbp} (2), where $\varphi: L_{2}\to L_{2}$ is the isomorphism of Lie algebras given by $x\mapsto x$, $y\mapsto y$, $z\mapsto z$ and $w\mapsto -x+w$.
\end{remark}

\begin{lemma}\label{nonilLb}
If a $p$-map on $L_{2}$ is not nilpotent then it
is  conjugate to one  of the following:
\begin{enumerate}
    \item $x^{[p]} = \alpha z,~ y^{[p]} = \beta_1 y +\beta_2 z, ~ z^{[p]} =0, ~ w^{[p]} = 0$, where $\alpha, \beta_1, \beta_2 \in \mathbb{F}$ with $\beta_1\neq 0$;
    \item $x^{[p]} = \alpha y, ~ y^{[p]} = \beta_1 y +\beta_2 z, ~ z^{[p]} =z, ~ w^{[p]} = 0$, where $\alpha, \beta_1, \beta_2 \in \mathbb{F}$.
\end{enumerate}
\end{lemma}

\begin{proof}
Let $[p]$ be a  $p$-map on $L_{2}$ that is not nilpotent. As explained in the paragraph before Proposition \ref{Lbp}, there exists $\alpha_i, \beta_i,\gamma_i, \delta_i\in \mathbb{F}$ for $i=1,2$ such that
\[
x^{[p]}=\alpha_{1}y+\alpha_{2}z, \quad y^{[p]}=\beta_{1}y+\beta_{2}z, \quad z^{[p]}=\gamma_{1}y+\gamma_{2}z, \quad w^{[p]}=\delta_{1}y+\delta_{2}z.
\]
Consider the isomorphism $\chi \in \mAut(L_{2})$ given by 
 \[
 x\mapsto x,\quad y\mapsto y, \quad z\mapsto z-uy,\quad w\mapsto w,
 \]
where $u\in \mathbb{F}$ is any scalar such that $\gamma_{1}+u^{p}\beta_{1}-u\gamma_{2}-u^{p+1}\beta_{2}=0.$
Such $u$ do exist since $\beta_{1},\beta_{2},\gamma_{2}$ cannot be $0$ simultaneously. Let $[p]_1:= \chi\circ [p] \circ \chi^{-1}$. A direct calculation shows that
\begin{align*}
  x^{[p]_1}= \alpha_1' y+ \alpha_2 z, \quad  y^{[p]_1}= \beta_1' y+\beta_2 z, \quad  z^{[p]_1} = \gamma z, \quad  w^{[p]_1}= \delta_1' y+\delta_2z,
\end{align*}
where $\alpha_1'= \alpha_1-u\alpha_2$, $\beta_1'=\beta_1 - u\beta_2$, $\gamma = \gamma_{2}+u^{p}\beta_{2}$, $\delta_1'= \delta_1-u\delta_2$. 
\begin{enumerate}
    \item Assume $\gamma=0$. Then $\beta_1'\neq 0$. Let $\varphi_1 \in \mAut(L_{2})$ be the isomorphism given by
 \[
 x\mapsto x+\left(\frac{\alpha_{1}'}{\beta_{1}'}\right)^{\frac{1}{p}}y,\quad y\mapsto y,\quad z\mapsto z,\quad w\mapsto w+\left(\frac{\delta_{1}'}{\beta_{1}'}\right)^{\frac{1}{p}}y.
 \]
 Then  $[p]_2:=\varphi_1\circ[p]_1\circ \varphi_1^{-1}$  satisfies that
 \[
 x^{[p]_2}=\alpha z,\quad y^{[p]_2}=\beta_{1}'y+\beta_{2}z,\quad z^{[p]_2}=0,\quad w^{[p]_2}=\delta z,
 \]
 where $\alpha= \alpha_{2}-\frac{\alpha_{1}'\beta_{2}}{\beta_{1}'}$ and $\delta=\delta_{2}-\frac{\delta_{1}'\beta_{2}}{\beta_{1}'}$.  We may assume $\delta\neq0$, or otherwise $[p]_2$ is already  of one of the desired forms. We proceed to split the discussion as following.
 \begin{itemize}
     \item Assume $\alpha=0$. Then $[p]_3:=\Delta\circ[p]_2\circ\Delta^{-1}$ satisfies that
\[
x^{[p]_3}=\delta z,\quad y^{[p]_3}=\beta_{1}'y+\beta_{2}z,\quad z^{[p]_3}=0,\quad w^{[p]_3}=\alpha z=0.
\] 
Here, $\Delta:L_{2}\to L_{2}$ is the isomorphism of Lie algebras given by
\[
x\mapsto -w,\quad y\mapsto -y,\quad z\mapsto -z,\quad w\mapsto -x.
\]
\item Assume $\alpha\neq 0$. For any $v_{1}, v_{2}\in \mathbb{F}$, define $\varphi_{2}\in \mAut(L_{2})$ by
\[
x\mapsto x,\quad y\mapsto y,\quad z\mapsto z,\quad w\mapsto -v_{1}x-v_{2}y+w.
\]
Then $[p]_3:= \varphi_{2}\circ[p]_2\circ\varphi_{2}^{-1}$ satisfies that 
\[
x^{[p]_3} =  \alpha z, \quad  y^{[p]_3} = \beta_{1}'y+\beta_{2}z, \quad z^{[p]_3}  = 0
\]
and
\[
w^{[p]_3}=
\begin{cases}
    (\delta+v_{1}^{p}\alpha+v_{2}^{p}\beta_{2})z+(v_{2}^{p}\beta_{1}'-v_{1})y, & \text{ when $p=2$}\\
    (\delta+v_{1}^{p}\alpha+v_{2}^{p}\beta_{2})z+v_{2}^{p}\beta_{1}'y,& \text{ when $p\geq3$}.
\end{cases}
\]
By choosing $v_{1}$ and $v_{2}$ appropriately one may get that $w^{[p]_3}=0$ as desired.
 \end{itemize}
 \item Now assume $\gamma\neq 0$. Let $\psi_1\in \mAut(L_{2})$ be the isomorphism given by
\[
x\mapsto x+\left(\gamma^{\frac{1}{p-1}}\alpha_{2}\right)^{\frac{1}{p}}z, \quad y\mapsto y, \quad z\mapsto\gamma^{\frac{1}{p-1}}z, \quad w\mapsto \left(\gamma^{\frac{1}{p-1}}\delta_{2}\right)^{\frac{1}{p}}z+w.
\]
Then  $[p]_2:=\psi_1\circ[p]_1\circ\psi_1^{-1}$ satisfies that
\[
x^{[p]_2}=\alpha_{1}' y,\quad y^{[p]_2}=\beta_{1}'y+ \left(\beta_{2}\gamma^{\frac{1}{p-1}}\right)z,\quad z^{[p]_2}=z,\quad w^{[p]_2}=\delta_{1}' y.
\]
For any $r\in \mathbb{F}$, let $\psi_{2}\in \mAut(L_{2})$ be the isomorphism given by
    \[
    x\mapsto x,\quad y\mapsto y,\quad z\mapsto z,\quad w\mapsto -rx+w.
    \]
    Then $[p]_3:= \psi_{2}\circ[p]_2\circ\psi_{2}^{-1}$ satisfies that  
    \[
    x^{[p]_3}=\alpha_{1}' y,\quad y^{[p]_3}=\beta_{1}'y+ \left(\beta_{2}\gamma^{\frac{1}{p-1}}\right)z,\quad z^{[p]_3}=z
    \]
and
    \[
    w^{[p]_3}=
    \begin{cases}
        (r^{p}\alpha_1'-r+\delta_1')y, &\text{when $p=2$},\\
        (r^{p}+\delta_1')y, &\text{when $p\geq3$}.
    \end{cases}
    \]
    By choosing  $r$ appropriately one may get that $w^{[p]_3}=0$ as desired.
\end{enumerate}
As a summary, in each case  we get the desired form up to conjugation of $p$-maps.
\end{proof}

\begin{lemma}\label{p7conjugate}
For any $\lambda\in \mathbb{F}$, let $[p]_{\lambda}$ be the $p$-map on $L_{2}$ that given in Proposition \ref{Lbp} (15). Then for any $\lambda, \mu \in \mathbb{F}$,  the $p$-maps $[p]_{\lambda}$ and $[p]_{\mu}$ are conjugate if and only if
\begin{eqnarray} \label{L50}
\lambda^{p(p-1)}(\mu^{p-1}+1)^{p+1} =\mu^{p(p-1)}(\lambda^{p-1}+1)^{p+1}.
\end{eqnarray}
\end{lemma}

\begin{proof}
Suppose $[p]_{\lambda}$ and $[p]_{\mu}$ are conjugate. Let $\varphi:L_{2}\rightarrow L_{2}$ be an isomorphism of Lie algebras such that $\varphi\circ [p]_{\lambda}= [p]_{\mu}\circ \varphi$. Describe $\varphi(x), \varphi(y), \varphi(z), \varphi(w)$ as that in \eqref{AutoLb} and let $b_{2}=a_{1}d_{4}-a_{4}d_{1}$. Then
considering $\varphi(y^{[p]_{\lambda}})=\varphi(y)^{[p]_{\mu}}$ and $\varphi(z^{[p]_{\lambda}})=\varphi(z)^{[p]_{\mu}}$, one has
\begin{eqnarray}\label{L51}
  c_{2}=c_{2}^{p}, \quad  \lambda c_{3}=b_{2}^{p}\mu, \quad  \lambda c_{2}=b_{2}^{p}-b_{2}, \quad \mu c_{2}^{p}=c_{3}-c_{3}^{p}. 
\end{eqnarray}
Note that $ b_2c_3-(b_2c_3)^p = b_2^p(c_3-c_3^p) - c_3(b_2^p-b_2)= b_2^p\mu c_2^p- c_3\lambda c_2 =0$. Thus
\begin{eqnarray}
 (b_2c_3)^{p-1} =1. \label{L55}
\end{eqnarray}
Clearly, $\lambda=0$ if and only if $\mu = 0$. Moreover, if  $c_{2}=0$, then  $b_{2},c_{3}\in \mathbb{F}_{p}^{\times}$ and hence $\lambda^{p-1}=\mu^{p-1}$. Therefore, to see the desired equality \eqref{L50}, we may assume 
\[
\lambda, ~ \mu, ~ c_2 \neq 0.
\]
Note that $c_2^{p^2-1}=c_2^{p-1}=1$. Then one has
\begin{eqnarray}
  &\left(\frac{\lambda}{\mu}\right)^{p-1}  \xlongequal{\eqref{L51}} \left(\frac{\lambda c_2}{\mu c_2^{p}}\right)^{p-1} \xlongequal{\eqref{L51}} \frac{b_{2}^{p-1}(b_{2}^{p-1}-1)^{p-1}}{c_{3}^{p-1}(1-c_{3}^{p-1})^{p-1}} \xlongequal{\eqref{L55}} b_{2}^{p^{2}-1},& \label{L56} \\
&(\lambda^{p}+\lambda)c_{2} \xlongequal{\eqref{L51}}  (\lambda c_2)^p +\lambda c_2  \xlongequal{\eqref{L51}}  b_2^{p^2} -b_2^p +\lambda c_2 \xlongequal{\eqref{L56}}\left (\frac{\lambda}{\mu}\right)^{p-1}b_{2}-b_{2},& \label{L57} \\
&(\lambda^{p}+\lambda)^{p^{2}-1}=\left( (\lambda^{p}+\lambda)c_{2}\right)^{p^2-1} \xlongequal{\eqref{L56}, \eqref{L57}} \left(\left(\frac{\lambda}{\mu}\right)^{p-1}-1 \right)^{p^{2}-1}\left(\frac{\lambda}{\mu}\right)^{p-1}.& \label{L58}
\end{eqnarray}
By \eqref{L58}, one readily has
\begin{equation}
    (\lambda^{p-1}+1)^{p^{2}-1}\lambda^{p(p-1)}\mu^{p^{2}(p-1)}=(\lambda^{p-1}-\mu^{p-1})^{p^{2}-1}. \label{L510}
\end{equation}
Since the symmetry of $\lambda$ and $\mu$, one also has
\begin{equation}
    (\mu^{p-1}+1)^{p^{2}-1}\mu^{p(p-1)}\lambda^{p^{2}(p-1)}=(\lambda^{p-1}-\mu^{p-1})^{p^{2}-1}. \label{L511}
\end{equation}
Calculate  $(\lambda^{p-1} +1)*\eqref{L510} - (\mu^{p-1}+1) * \eqref{L511}$, one obtains
\begin{equation}
    (\lambda^{p-1}+1)^{p^{2}}\lambda^{p(p-1)}\mu^{p^{2}(p-1)}-(\mu^{p-1}+1)^{p^{2}}\mu^{p(p-1)}\lambda^{p^{2}(p-1)}=(\lambda^{p-1}-\mu^{p-1})^{p^{2}}.
\end{equation}
Then take the $p$-th root, it follows that
\begin{equation}
(\lambda^{p-1}+1)^{p}\lambda^{p-1}\mu^{p(p-1)}-(\mu^{p-1}+1)^{p}\lambda^{p(p-1)}\mu^{(p-1)}=(\lambda^{p-1}-\mu^{p-1})^{p}.
\end{equation}
Consequently,
\begin{equation}
\lambda^{p(p-1)}(\mu^{p^{2}-1}+\mu^{p-1}+1) =\mu^{p(p-1)}(\lambda^{p^{2}-1}+\lambda^{p-1}+1).
\end{equation}
The equality  \eqref{L50} follows, since $(\lambda^{p-1}+1)^{p+1} = \lambda^{p^{2}-1}+\lambda^{p(p-1)} + \lambda^{p-1}+1$ and similar for $\mu$.

Next we show the converse implication. Suppose the equality \eqref{L50} holds for $\lambda, \mu \in \mathbb{F}$. Clearly, to see that $[p]_{\lambda}$ and $[p]_{\mu}$ are conjugate, we may assume $\lambda, \mu\neq 0$. Note that for any scalars $b_2, c_2, c_3\in \mathbb{F}$ with $b_2,c_3\neq 0$, the linear map  $\psi:L_{2}\to L_{2}$ given by
\[
x\mapsto b_2x, \quad y\mapsto b_2 y, \quad z\mapsto c_2y +c_3 z, \quad w\mapsto w
\]
is an isomorphism of Lie algebras. It remains to specify $b_2, c_2, c_3$ such that $\psi\circ [p]_{\lambda}= [p]_{\mu}\circ \psi$. It is easy to check  that $\psi(x^{[p]_{\lambda}}) =\psi(x)^{[p]_{\mu}}=0$ and $\psi(w^{[p]_{\lambda}}) =\psi(w)^{[p]_{\mu}} =0$, and if the equalities in \eqref{L51} hold then $\psi(y^{[p]_{\lambda}}) =\psi(y)^{[p]_{\mu}}$ and $\psi(z^{[p]_{\lambda}}) =\psi(z)^{[p]_{\mu}}$,  for any $b_2, c_2, c_3$. Therefore, it suffices to specify $b_2, c_2, c_3$ with $b_2,c_3\neq 0$ such that the equalities in \eqref{L51} hold. If $\lambda^{p-1}=-1$ then $\mu^{p-1}=-1$ by the equality \eqref{L50}. Clearly, in this case we may take
\[
b_2=1, \quad c_2=0, \quad c_3= \frac{\mu}{\lambda}.
\]
Otherwise, if $\lambda^{p-1}\neq -1$ then we may take
\[
b_{2}=\left(\frac{\lambda}{\mu}\right)^{\frac{1}{p+1}}\xi,\quad  c_{2}=\frac{\lambda^{\frac{p-1}{p+1}}\xi^{p}-\mu^{\frac{p-1}{p+1}}\xi}{(\lambda\mu)^{\frac{p}{p+1}}},\quad c_{3}=\left(\frac{\mu}{\lambda}\right)^{\frac{1}{p+1}}\xi^{p}, 
\]
where $\xi\in \mathbb{F}$ is any scalar satisfying $\xi^{p-1}=\left(\frac{\lambda}{\mu}\right)^{\frac{p(p-1)}{p+1}}\left(\frac{\mu^{p-1}+1}{\lambda^{p-1}+1}\right)$. Indeed, the equality \eqref{L50} implies that $\xi^{p^{2}}= (\xi^{p-1})^{p+1} \xi=\xi$. Consequently, one has
\begin{eqnarray*}
c_{2} - c_{2}^{p} &=& \frac{\left(\frac{\lambda}{\mu}\right)^{\frac{p}{p+1}}\xi^{p}-\left(\frac{\lambda}{\mu}\right)^{\frac{1}{p+1}}\xi}{\lambda} - \frac{\left(\frac{\lambda}{\mu}\right)^{\frac{p^{2}}{p+1}}\xi^{p^{2}}-\left(\frac{\lambda}{\mu}\right)^{\frac{p}{p+1}}\xi^{p}}{\lambda^{p}}  \\
 &=& \frac{\xi}{\lambda^p}\left( \left(\frac{\lambda}{\mu}\right)^{\frac{p}{p+1}}(\lambda^{p-1}+1)\xi^{p-1} - \left(\frac{\lambda}{\mu}\right)^{\frac{p^{2}}{p+1}}(\mu^{p-1}+1) \right) \\
 &= &0
\end{eqnarray*}
The remaining three equalities in \eqref{L51} are easy to verify.    
\end{proof}

\begin{proof}[Proof of Proposition \ref{Lbp}]
By \cite[Theorem 2.1]{schneider2016classification}, the $p$-maps in Parts (1) to (8) form the conjugacy classes of the nilpotent $p$-maps on $L_{2}$. We proceed to show the remaining ones form the conjugacy classes of the  $p$-maps on $L_{2}$ that are not nilpotent. To this end, we denote by $[p]_{9}, \ldots, [p]_{14}$ and $[p]_{15,\lambda}$, where $\lambda\in \mathbb{F}$, the $p$-maps listed in Parts (9) to (15) correspondingly. We ignore the restriction on $p$ in Parts (11) and (13) yet. Clearly, they all are not nilpotent.

Let $[p]$ be any $p$-map on $L_{2}$ that is not nilpotent. According to Lemma \ref{nonilLb}, we may assume $[p]$ satisfies one of the following two conditions: (I)  $x^{[p]} = \alpha z$, $y^{[p]} = \beta_1 y +\beta_2 z$, $z^{[p]} =0$ and $w^{[p]} = 0$ for some $\alpha,  \beta_1, \beta_2 \in \mathbb{F}$ with $\beta_1\neq 0$; (II) $x^{[p]} = \alpha y$, $y^{[p]} = \beta_1 y +\beta_2 z$, $z^{[p]} =z$ and $w^{[p]} = 0$ for some $\alpha,  \beta_1, \beta_2 \in \mathbb{F}$.  We proceed to split the discussion as following.
\begin{enumerate}
\item Assume (I) holds and $\beta_2=\alpha=0$. Then $\sigma_{a,b,1,1}\circ[p]\circ\sigma_{a,b,1,1}^{-1}=[p]_{9}$, where $a=b=\beta_{1}^{\frac{1}{p-1}}$.
\item Assume (I) holds,  $\beta_2=0$ and $\alpha\neq 0$. Then $\sigma_{a,b,1,d}\circ[p]\circ\sigma_{a,b,1,d}^{-1}=[p]_{10}$, where $a=\alpha^{\frac{1}{p}}$, $b=\beta_{1}^{\frac{1}{p-1}}$, $d=\beta_{1}^{\frac{1}{p-1}}\alpha^{-\frac{1}{p}}$.
\item Assume (I) holds, $\beta_2\neq0$ and $\alpha=0$. Then $\sigma_{a,b,c,1}\circ[p]\circ\sigma_{a,b,c,1}^{-1}=[p]_{11}$, where $a=b=\beta_{1}^{\frac{1}{p-1}}$, $c=\beta_{2}^{-1}\beta_{1}^{\frac{p}{p-1}}$. 
\item Assume (I) holds, $\beta_2\neq0$ and $\alpha\neq0$. Then $\sigma_{a,b,c,d}\circ[p]\circ\sigma_{a,b,c,d}^{-1}=[p]_{12}$, where $a=\left(\frac{\alpha}{\beta_{2}}\right)^{\frac{1}{p}}\beta_{1}^{\frac{1}{p-1}}$, $b=\beta_{1}^{\frac{1}{p-1}}$, $c=\beta_{1}^{\frac{p}{p-1}}\beta_{2}^{-1}$, $d=\left(\frac{\beta_{2}}{\alpha}\right)^{\frac{1}{p}}$.
\item Assume (II) holds and $\beta_1=\beta_2=\alpha=0$. Then $[p]=[p]_{13}$.
\item Assume (II) holds, $\beta_1=\beta_2=0$ and $\alpha\neq 0$. Then  $\sigma_{a,1,1,d}\circ[p]\circ\sigma_{a,1,1,d}^{-1} =[p]_{14}$, where $a=\alpha^{\frac{1}{p}}$ and $d=\alpha^{-\frac{1}{p}}$.
\item Assume (II) holds, $\beta_1=0$ and $\beta_2\neq 0$. Then $[p]':=\tau_{1}\circ[p]\circ\tau_{1}^{-1}$ satisfies that 
\[
x^{[p]'} = -\alpha \beta_2^{\frac{1}{p}} z, \quad y^{[p]'} = \beta_2^{\frac{p-1}{p}} y +\beta_2 z, \quad z^{[p]'} =0, \quad w^{[p]'}=0.
\]
Here, $\tau_{1}:L_{2}\to L_{2}$ is the isomorphism of Lie algebras given by
\[
x\mapsto x+ \left(\alpha^{\frac{1}{p}} \beta_{2}^{-\frac{p-1}{p^2}}\right) y,\quad y\mapsto y,\quad z\mapsto z+\beta_{2}^{-\frac{1}{p}}y,\quad w\mapsto w.
\]
Clearly, $[p]'$ falls into the situations (3) or (4).
\item Assume (II) holds, $\beta_1\neq0$ and $\alpha=0$. Then $\sigma_{a,b,1,1}\circ[p]\circ\sigma_{a,b,1,1}^{-1}
=[p]_{15,\lambda'}$ where $a=b=\beta_{1}^{\frac{1}{p-1}}$ and $\lambda'=\beta_{1}^{\frac{p}{1-p}}\beta_{2}$.
\item Assume (II) holds, $\beta_1\neq 0$ and $\alpha\neq 0$. Then $[p]':= \tau_2\circ [p] \circ \tau_2^{-1}$ satisfies that 
\[
x^{[p]'}=0,\quad y^{[p]'}=\beta_{1}y+\beta_{2}z,\quad z^{[p]'}=z,\quad w^{[p]'}=0.
\]
Here, $\tau_{2}:L_{2}\to L_{2}$ is the isomorphism of Lie algerbas given by
\[
    x\mapsto \alpha^{\frac{1}{p}}x+ \left(\frac{\alpha}{\beta_{1}}\right)^{\frac{1}{p}}y-\left(\frac{\alpha\beta_{2}}{\beta_{1}}\right)^{\frac{1}{p}}z,\quad y\mapsto y,\quad
    z\mapsto z,\quad w\mapsto \alpha^{-\frac{1}{p}}w.
\]
Clearly, $[p]'$ falls into the situation (8).
\end{enumerate}
In addition, if $p=2$ then $\varphi\circ[p]_{11}\circ\varphi^{-1}= [p]_{12}$ and $\psi\circ[p]_{13}\circ \psi^{-1}=[p]_{14}$, where $\varphi: L_{2}\to L_{2}$ (resp. $\psi: L_{2}\to L_{2}$) is the isomorphism of Lie algebras given by 
 \[
x\mapsto x-y-w  ~~(\text{resp. }  x-w),\quad y\mapsto y,\quad z\mapsto z,\quad w\mapsto w.
\]
Here, we used the equality $(x+w)^{[2]}=x^{[2]}+w^{[2]}+y$ for any $2$-map $[2]$ on $L_{2}$. As a summary, we showed that every non-nilpotent $p$-map on $L_{2}$ is conjugate to one in Parts (9) to (15).

Next, we  show they are not conjugate to each other except those indicated in Part (15). By a straightforward but tedious calculation on $\dim L_{2}^{[p]}$, $\dim Z(L_{2})^{[p]}$  and $\dim [L_{2},L_{2}]^{[p]}$, we obtain Table \ref{Lbpdim}.
\begin{table}[htbp]
\centering
\caption{3-tuple of invariant of $L=L_{2}$}\label{Lbpdim}
\setlength{\tabcolsep}{3.5mm}{
\renewcommand\arraystretch{1.4}
\begin{tabular}{c|ccc|c}
\hline
$[p]$      & $\dim L^{[p]}$ & $\dim Z(L)^{[p]}$ &$\dim [L,L]^{[p]}$  & Condition  \\ \hline
$[p]_{9}$  & 1 & 1 & 1 & \\\hline
$[p]_{10}$ & 2 & 1 & 1 &            \\ \hline
$[p]_{11}$  & 1 & 1 & 1 & $p\geq3$     \\ \hline
$[p]_{12}$ & 2 & 1 & 1 &  \\\hline
$[p]_{13}$  & 1 & 1 & 0 & $p\geq3$  \\ \hline
$[p]_{14}$ & 2 & 1 & 0 &  \\\hline
$[p]_{15,\lambda}$ & 2 & 2 & 1 &      \\    \hline  
\end{tabular}}
\end{table}
So by Lemma \ref{ConjInvar} and Lemma \ref{p7conjugate}, 
it remains to show  $[p]_{9}$ and $[p]_{11}$ as well as $[p]_{10}$ and $[p]_{12}$ are not conjugate. 
But if there exists an isomorphism $\Phi\in \mAut(L_{2})$ such that $\Phi\circ [p]_{9}\circ \Phi^{-1} = [p]_{11}$ or $\Phi\circ [p]_{10}\circ \Phi^{-1} = [p]_{12}$, then taking value at $y$ one may get in each case that $$b_{2}^{1-p} y =y+z,$$ which is impossible. Here, $b_{2}\in \mathbb{F}$ is defined by $\Phi(y) =b_{2}y$.
\end{proof}

\subsection{Case $L_{3}$}

Recall that the Lie algebra $L_{3}$ has a basis $x, y, z, w$ with structure relations $[w,x]=y$ and $[w,y]=z$.  Clearly, $Z(L_{3}) =\mspan\{z\}$. In addition, a linear map $\varphi: L_{3} \to L_{3}$ is a homomorphism of Lie algebras if and only if  $\varphi(z) \in \mspan\{z\}$, $[\varphi(w), \varphi(x)] = \varphi(y)$ and $[\varphi(w), \varphi(y)] = \varphi(z)$. It follows that $\varphi \in \mAut(L_{3})$ if and only if it is given by 
\begin{equation}\label{AutoLc}
    \begin{aligned}
    x&\mapsto a_{1}x+a_{2}y+a_{3}z+a_{4}w, \quad  \quad \quad&
    y&\mapsto (a_{1}d_{4}-a_{4}d_{1})y + (a_{2}d_{4}-a_{4}d_{2})z,\\
    z&\mapsto (a_{1}d_{4}-a_{4}d_{1})d_{4}z, \quad \quad \quad & 
    w&\mapsto d_{1}x+d_{2}y+d_{3}z+d_{4}w
    \end{aligned}
    \end{equation}
for some $a_i, b_i, d_i\in \mathbb{F}$, which satisfy $(a_{1}d_{4}-a_{4}d_{1})d_{4}  \neq 0$.

\begin{prop}\label{Lcp2}
Let $p=2$. There is no $2$-map on $L_{3}$.
\end{prop}

\begin{proof}
Clearly, $(\mad w)^{2}(x)=z$. Thus $(\mad w)^{2}\neq \mad f$ for any $f\in L_{3}$, because 
$[f,x] \in \mspan\{ y\}$. The result follows directly.   
\end{proof}

Assume $p\geq 3$. It is easy to verify that $(\mad x)^{p} = (\mad y)^{p} = (\mad z)^{p} = (\mad w)^{p} =0$. Thus by Lemma \ref{pmapconstruct}, the $p$-maps on $L_{3}$ are the $p$-semilinear maps of the form:
\begin{equation}
x^{[p]}=\alpha z, \quad  y^{[p]}=\beta z, \quad z^{[p]}=\gamma z, \quad w^{[p]}=\delta z,
\end{equation}
where $\alpha, \beta, \gamma, \delta \in \mathbb{F}$. The following result classifies the $p$-maps on $L_{3}$ for $p\geq 3$. 

\begin{prop}\label{Lcp}
Let $p\geq 3$. The conjugacy classes of the $p$-maps on $L_{3}$ are as follows:   
\begin{enumerate}
    \item Trivial $p$-map, for $p\geq 5$;
    \item $x^{[p]}=0,~y^{[p]}=0,~z^{[p]}=0,~w^{[p]}=z$;
    \item $x^{[p]}=0,~y^{[p]}=z,~z^{[p]}=0,~w^{[p]}=0$;
    \item $x^{[p]}=z,~y^{[p]}=0,~z^{[p]}=0,~w^{[p]}=0$;
    \item $x^{[p]}=0,~y^{[p]}=0,~z^{[p]}=z,~w^{[p]}=0$.
\end{enumerate}
The $p$-maps in (1) to (4) are nilpotent, while the last one is not.
\end{prop}
\begin{proof}
By \cite[Theorem 2.1]{schneider2016classification}, the $p$-maps in Parts (1) to (4) form conjugacy classes of the nilpotent $p$-maps on $L_{3}$. It remains to show every non-nilpotent $p$-map $[p]$ on $L_{3}$ is conjugate to the one in Part (5). As explained above, we may assume $x^{[p]}=\alpha z$, $y^{[p]}=\beta z$, $z^{[p]}=\gamma z$ and $w^{[p]}=\delta z$ for some $\alpha, \beta, \gamma, \delta \in \mathbb{F}$.
Since $[p]$ is not nilpotent, $\gamma\neq 0$. Define $\varphi \in \mAut(L_{3})$ by
\[
x\mapsto s^{p} x+\beta^{\frac{1}{p}}sy+\alpha^{\frac{1}{p}}sz, \quad y\mapsto s^{p} y+\beta^{\frac{1}{p}}sz, \quad  z\mapsto s^{p}z, \quad w\mapsto \delta^{\frac{1}{p}}sz+w,
\]
where $s=\gamma^{\frac{1}{p(p-1)}}$.  Then $\varphi\circ [p]\circ \varphi^{-1}$ is exactly the $p$-map given in Part (5).
\end{proof}

\begin{remark}
Assume $p=3$. Nilpotent $3$-maps on $L_{3}$ are classified into $4$ classes in  \cite[Theorem 2.1 (4/3)]{schneider2016classification}, while they are only classified into $3$ classes in  Proposition \ref{Lcp}. The difference is due to that the base field in \cite{schneider2016classification} is only assumed to be perfect. In fact, the $3$-maps in $(c)$ and  $(d)$ of \cite[Theorem 2.1 (4/3)]{schneider2016classification} are conjugate to each other when the base field is algebraically closed. 
\end{remark}

\begin{remark}
Let $[p]$ be the $p$-map on $L_{3}$ given by $x^{[p]} =0$, $y^{[p]} =0$, $z^{[p]} =0 $ and $w^{[p]}=0$ for $p \geq 3$. It is the trivial $p$-map when $p\geq 5$. However, it is not for $p=3$. Indeed, one has 
\[
(x+w)^{[3]} =x^{[3]} +w^{[3]} + [w, [w,x]] +[x,[x,w]] = z.
\]
Actually, $\varphi\circ [3]\circ \varphi^{-1}$ is the one given in Proposition \ref{Lcp} (4), where $\varphi: L_{3}\to L_{3}$ is the isomorphism of Lie algebras given by $x\mapsto x$, $y\mapsto y$, $z\mapsto z$ and $w\mapsto -x+w$.
\end{remark}

\subsection{Case $L_{4}$}

Recall that the Lie algebra $L_{4}$ has a basis $x, y, z, w$ with structure relations $[w,x] =x$. Clearly, $Z(L_{4}) =\mspan\{y,z\}$. In addition, a linear map $\varphi: L_{4} \to L_{4}$ is a homomorphism of Lie algebras if and only if  $\varphi(y), \varphi(z) \in \mspan\{y, z\}$ and $[\varphi(w), \varphi(x)] = \varphi(x)$. It follows that $\varphi \in \mAut(L_{4})$ if and only if it is given by 
\begin{eqnarray}\label{AutoLd} 
x\mapsto a_{1}x, \quad y\mapsto b_{2}y + b_{3}z, \quad z\mapsto c_{2}y+c_{3}z, \quad w\mapsto d_{1}x+d_{2}y+d_{3}z+w  
\end{eqnarray}
for some $a_1, b_i, c_i, d_i\in \mathbb{F}$, which satisfy $a_{1}(b_{2}c_{3} -b_{3}c_{2}) \neq 0$.

It is easy to verify that
$(\mad x)^p = (\mad y)^p = (\mad z)^p =0$ and $(\mad w)^p = \mad w$. Thus by Lemma \ref{pmapconstruct}, the $p$-maps on $L_{4}$ are the $p$-semilinear maps of the form:
\begin{equation}
x^{[p]}=\alpha_{1}y+\alpha_{2}z, \quad y^{[p]}=\beta_{1}y+\beta_{2}z, \quad z^{[p]}=\gamma_{1}y+\gamma_{2}z, \quad w^{[p]}=\delta_{1}y+\delta_{2}z +w,
\end{equation}
where $\alpha_i, \beta_i, \gamma_i, \delta_i\in \mathbb{F}$. The next result classifies the $p$-maps on $L_{4}$.

\begin{prop}\label{Ldp}
The conjugacy classes of the $p$-maps on $L_{4}$ are as follows:
    \begin{enumerate}
    \item $x^{[p]}=0,~y^{[p]}=0,~z^{[p]}=0,~w^{[p]}=w$;
    \item $x^{[p]}=y,~y^{[p]}=0,~z^{[p]}=0,~w^{[p]}=w$;
    \item $x^{[p]}=0,~y^{[p]}=z,~z^{[p]}=0,~w^{[p]}=w$;
    \item $x^{[p]}=z,~ y^{[p]}=z,~z^{[p]}=0,~w^{[p]}=w$;
    \item $x^{[p]}=y,~ y^{[p]}=z,~z^{[p]}=0,~w^{[p]}=w$;
    \item $x^{[p]}=0,~y^{[p]}=y,~z^{[p]}=0,~w^{[p]}=w$;
    \item $x^{[p]}=z,~ y^{[p]}=y,~z^{[p]}=0,~w^{[p]}=w$;
    \item $x^{[p]}=y,~ y^{[p]}=y,~z^{[p]}=0,~w^{[p]}=w$;
    \item $x^{[p]}=y+z,~ y^{[p]}=y,~z^{[p]}=0,~w^{[p]}=w$;
    \item $x^{[p]}=0,~ y^{[p]}=y,~ z^{[p]}=z,~ w^{[p]}=w$;
    \item $x^{[p]}=y+ \lambda z,~ y^{[p]}=y,~ z^{[p]}=z,~ w^{[p]}=w$, $\lambda\in \mathbb{F}$. Here, $\lambda_1, \lambda_2\in \mathbb{F}$ represent conjugate $p$-maps if and only if  $\left(\begin{smallmatrix} 1 \\ \lambda_{1} \end{smallmatrix}\right)=aA\!\left(\begin{smallmatrix} 1 \\ \lambda_{2} \end{smallmatrix}\right)$ and $A\in \mathrm{GL}_{2}(\mathbb{F}_{p})$.
    \end{enumerate}
\end{prop}
The proof of this proposition will be addressed at the end of this subsection. 

\begin{lemma}\label{Ldconjugate}
Every $p$-map on $L_{4}$ is conjugate to one of the following:
\begin{enumerate}
\item $x^{[p]} =\lambda_1 y+\lambda_2 z, ~ y^{[p]}=0, ~z^{[p]}=0, ~ w^{[p]} = w$, where $\lambda_1,\lambda_2\in \mathbb{F}$;
    \item $x^{[p]} =\lambda_1 y+\lambda_2 z, ~ y^{[p]}=z, ~z^{[p]}=0, ~ w^{[p]} = w$, where $\lambda_1,\lambda_2\in \mathbb{F}$;
    \item $x^{[p]} =\lambda_1 y+\lambda_2 z, ~ y^{[p]}=y,~z^{[p]}=0, ~ w^{[p]} = w$, where $\lambda_1,\lambda_2\in \mathbb{F}$;
    \item $x^{[p]} =\lambda_1 y+\lambda_2 z, ~ y^{[p]}=y,~z^{[p]}=z, ~ w^{[p]} = w$, where $\lambda_1,\lambda_2\in \mathbb{F}$.
\end{enumerate}
\end{lemma}

\begin{proof}
Let $[p]$ be an arbitrary $p$-map on $L_{4}$. As explained in the paragraph before Proposition \ref{Ldp}, we may assume $x^{[p]}=\alpha_{1}y+\alpha_{2}z,$ $y^{[p]}=\beta_{1}y+\beta_{2}z,$ $z^{[p]}=\gamma_{1}y+\gamma_{2}z,$ and $w^{[p]}=\delta_{1}y+\delta_{2}z +w$ for some $\alpha_i, \beta_i, \gamma_i, \delta_i\in \mathbb{F}$.
Choose scalars $r,s\in \mathbb{F}$ so that 
\[
r=\delta_{1}+r^{p}\beta_{1}+s^{p}\gamma_{1} \quad \text{and} \quad s=\delta_{2}+r^{p}\beta_{2}+s^{p}\gamma_{2}.
\]
Clearly, such $r, s$ do exist. Let $\varphi \in \mAut(L_{4})$ is the isomorphism given by
\[
x\mapsto x,\quad y\mapsto y,\quad z\mapsto z,\quad w\mapsto w+ry+sz.
\]
Then for $[p]_1:= \varphi\circ [p] \circ \varphi^{-1}$, one has $x^{[p]_1} =x^{[p]}$,  $y^{[p]_1} =y^{[p]}$, $z^{[p]_1} =z^{[p]}$ and
\[
w^{[p]_1} = \varphi\left((w-ry-sz)^{[p]}\right) = \varphi\left(w^{[p]} - r^{p} y^{[p]} - s^{p} z^{[p]}\right) =w.
\]
For any $A\in \mathrm{GL}_{2}(\mathbb{F})$, define $\psi_A\in \mAut(L_{4})$  by 
\[
x\mapsto x, \quad (y, z) \mapsto (y, z)A, \quad w\mapsto w.
\]
Let $[p]_A:=\psi_A\circ [p]_1\circ \psi_A^{-1}$.
For any $A=\begin{pmatrix}[1.2]
        a_{11} & a_{12}  \\
        a_{21} & a_{22}  
    \end{pmatrix}\in \mathrm{GL}_{2}(\mathbb{F})$, one has 
\[
x^{[p]_A}= (y,z) A \begin{pmatrix}[1.2]
        \alpha_{1}  \\
        \alpha_{2}   
    \end{pmatrix}, \quad (y^{[p]_A}, z^{[p]_A}) = (y,z)A \begin{pmatrix}[1.2]
        \beta_{1} & \gamma_{1}  \\
        \beta_{2} & \gamma_{2}  
    \end{pmatrix} \tilde{A}, \quad w^{[p]_A} =w.
\]
where  $\tilde{A}= \frac{1}{\det(A)^{p}} \begin{pmatrix}[1.2]
        a_{22}^p & -a_{12}^p  \\
        -a_{21}^p & a_{11}^p  
    \end{pmatrix}$. By the classification result of the $p$-maps on $\mspan\{y, z\}$, which is considered as the abelian Lie algebra of dimension $2$, there exists $A\in \mathrm{GL}_{2}(\mathbb{F})$ such that
\begin{eqnarray*}
 A \begin{pmatrix}[1.2]
        \beta_{1} & \gamma_{1}  \\
        \beta_{2} & \gamma_{2}  
    \end{pmatrix} \tilde{A} &\in& \big \{\begin{pmatrix}[1.2]
        0 & 0  \\
        0 & 0  
    \end{pmatrix}, ~ \begin{pmatrix}[1.2]
        0 & 0  \\
        1 & 0  
    \end{pmatrix}, ~ \begin{pmatrix}[1.2]
        1 & 0  \\
        0 & 0  
    \end{pmatrix}, ~ \begin{pmatrix}[1.2]
        1 & 0  \\
        0 & 1  
    \end{pmatrix} \big \}.   
\end{eqnarray*}
See \cite[Proposition A.1]{wang2015classification} for a reference. The result then follows immediately.
\end{proof}

\begin{lemma}\label{p10conjugate}
For any $\lambda =\left(\begin{smallmatrix} \lambda_{1} \\ \lambda_{2} \end{smallmatrix}\right) \in \mathbb{F}\times \mathbb{F}$, let $[p]_{\lambda}$ be the $p$-map on $L_{4}$ that given in Part (4) of Lemma \ref{Ldconjugate}. 
Then for any $\lambda, \mu \in \mathbb{F}\times \mathbb{F}$, the $p$-maps $[p]_{\lambda}$ and $[p]_{\mu}$ are conjugate if and only if $\lambda = aA\mu$ for some nonzero element $a\in \mathbb{F}^\times$ and invertible matrix $A\in \mathrm{GL}_{2}(\mathbb{F}_{p})$.
\end{lemma}

\begin{proof}
Firstly, we show the forward implication. Suppose $[p]_{\lambda}$ and $[p]_{\mu}$ are conjugate. Let $\varphi\in \mAut(L_{4})$ be an isomorphism such that $\varphi\circ [p]_{\lambda} \circ \varphi^{-1} = [p]_{\mu}$. Describe $\varphi(x), \varphi(y), \varphi(z), \varphi(w)$  as that in \eqref{AutoLd}. Consider $\varphi(x)^{[p]_{\lambda}}=\varphi(x^{[p]_{\mu}})$, one obtains 
\[
\begin{pmatrix}[1.2]
        \lambda_1  \\
        \lambda_2   
    \end{pmatrix}
    = a_{1}^{-p}  \begin{pmatrix}[1.2]
        b_{2} & c_{2}  \\
        b_{3} & c_{3}  
    \end{pmatrix} \begin{pmatrix}[1.2]
        \mu_1  \\
        \mu_2   
    \end{pmatrix}.
\]
In addition, 
consider $\varphi(y)^{[p]_{\lambda}}=\varphi(y^{[p]_{\mu}})$ and $\varphi(z)^{[p]_{\lambda}}=\varphi(z^{[p]_{\mu}})$, one has
\[
b_{2}^{p}=b_{2}, \quad b_{3}^{p}=b_{3}, \quad c_{2}^{p}=c_{2}, \quad c_{3}^{p}=c_{3}. 
\]
So $b_{2},b_{3},c_{2},c_{3}\in\mathbb{F}_{p}$ and hence the desired condition holds.

Conversely, assume $\lambda = a A\mu$ for some nonzero element $a\in \mathbb{F}^\times$ and invertible matrix $A\in \mathrm{GL}_{2}(\mathbb{F}_{p})$. Let $\varphi:L_{4}\to L_{4}$ be the linear map given by $x\mapsto a^{-\frac{1}{p}}x$, $(y,z)\mapsto (y,z)A$ and $w\mapsto w$. Readily, it is an isomorphism of Lie algebras and $\varphi\circ [p]_{\lambda} \circ \varphi^{-1} = [p]_{\mu}$. 
\end{proof}

\begin{proof}[Proof of Proposition \ref{Ldp}]
Let $[p]$ be a $p$-map on $L_{4}$. By Lemma \ref{Ldconjugate}, we may assume $x^{[p]} = \lambda_{1}y+ \lambda_{2} z$ for some $\lambda_{1}, \lambda_{2}\in \mathbb{F}$, $w^{[p]}=w$ and  one of the following conditions hold: (I) $y^{[p]}=0,~z^{[p]}=0$; (II) $y^{[p]}=z,~z^{[p]}=0$; (III) $y^{[p]}=y,~z^{[p]}=0$; (IV) $y^{[p]}=y,~z^{[p]}=z$. Let $[p]_{1}, \ldots, [p]_{10}$ and $[p]_{11,\lambda}$, $\lambda\in \mathbb{F}$, be the $p$-maps listed in Proposition \ref{Ldp} correspondingly.
\begin{enumerate}
\item Assume (I) hold. If $\lambda_{1}=\lambda_{2}=0$ then  $[p]=[p]_{1}$; If $\lambda_{1}=0$ and $\lambda_{2}\neq0$ then $\varphi\circ [p]\circ \varphi^{-1} =[p]_{2}$, where $\varphi\in \mAut(L_{4})$ is the isomorphism given by 
\[
x\mapsto \lambda_{2}^{\frac{1}{p}}x,\quad y\mapsto z, \quad z\mapsto y, \quad w\mapsto w;
\]
If $\lambda_{1}\neq 0$ then $\psi\circ [p] \circ \psi^{-1}= [p]_{2}$, where $\psi\in \mAut(L_{4})$ is the isomorphism given by 
\[
x\mapsto x,\quad y\mapsto \frac{1}{\alpha_{1}}y-\frac{1}{\alpha_{1}^{p+1}}\alpha_{2}z,\quad z\mapsto \frac{1}{\alpha_{1}^{p}}z,\quad w\mapsto w.
\]
\item Assume (II) hold. If $\lambda_{1}=\lambda_{2}=0$ then $[p]= [p]_{3}$; If $\lambda_{1}=0$ and $\lambda_{2}\neq0$ then $\sigma_{a,1,1,1}\circ[p]\circ\sigma_{a,1,1,1}^{-1} = [p]_{4}$, where $a=\lambda_{2}^{\frac{1}{p}}$; If $\lambda_{1}\neq 0$ then $\psi\circ [p] \circ \psi^{-1}= [p]_{5}$.
\item Assume (III) hold. If $\lambda_{1}=\lambda_{2}=0$ then $[p]= [p]_{6}$; If $\lambda_{1}=0$ and $\lambda_{2}\neq0$ then $\sigma_{a,1,1,1}\circ[p]\circ\sigma_{a,1,1,1}^{-1} = [p]_{7}$, where $a=\lambda_{2}^{\frac{1}{p}}$; If $\lambda_{1}\neq0$ and $\lambda_{2}=0$ then $\sigma_{a,1,1,1}\circ[p]\circ\sigma_{a,1,1,1}^{-1} = [p]_{8}$, where $a=\lambda_{1}^{\frac{1}{p}}$;
If $\lambda_{1}, \lambda_{2} \neq0$ then $\sigma_{a,1,c,1}\circ [p]\circ\sigma_{a,1,c,1}^{-1}=[p]_{9}$, where $a=\lambda_{1}^{\frac{1}{p}}$ and $c=\frac{\lambda_{1}}{\lambda_{2}}$.
\item Assume (IV) hold. If $\lambda_{1}=\lambda_{2}=0$ then $[p]= [p]_{10}$. Otherwise, if $\lambda_{1}\neq0$ or $\lambda_{2}\neq0$,  then it is easy to see that there exists $\lambda\in \mathbb{F}$, $a \in \mathbb{F}^\times$ and $A\in \mathrm{GL}_2(\mathbb{F}_{p})$ such that 
$(1,\lambda)^{t} =a A(\lambda_{1}, \lambda_{2})^t$, and therefore $[p]$ is conjugate to $[p]_{11,\lambda}$ by Lemma \ref{p10conjugate}.
\end{enumerate}
As a summary, we showed that every $p$-map on $L_{4}$ is conjugate to one of the listed $p$-maps. By a straightforward but tedious calculation on $\dim L_{4}^{[p]}$, $\dim[L_{4},L_{4}]^{[p]}$, $\dim Z(L_{4})^{p}$, $\dim Z(L_{4})^{[p]^{2}}$ and $\dim[L_{4},L_{4}]^{[p]^{2}}$, we have Table \ref{Ldpdim}. 
\begin{table}[htbp]
\centering
\caption{5-tuple of invariant of $L=L_{4}$}\label{Ldpdim}
\setlength{\tabcolsep}{3.3mm}{
\renewcommand\arraystretch{1.35}
\begin{tabular}{cccccc}
\hline
$[p]$      & $\dim L^{[p]}$ & $\dim [L,L]^{[p]}$ & $\dim Z(L)^{[p]}$ & $\dim Z(L)^{[p]^{2}}$ &$\dim [L,L]^{[p]^{2}}$ \\ \hline
$[p]_{1}$  & 2 & 0 & 0 & 0 & 0 \\\hline
$[p]_{2}$  & 3 & 1 & 0 & 0  & 0\\ \hline
$[p]_{3}$  & 3 & 0 &1  & 0  &0  \\ \hline
$[p]_{4}$ & 3 & 1 & 1 & 0   & 0 \\\hline
$[p]_{5}$ & 4 & 1 & 1 &  0   & 0       \\ \hline
$[p]_{6}$  & 3 & 0 &1  & 1  & 0 \\ \hline 
$[p]_{7}$ & 4 & 1 & 1 & 1    & 0      \\ \hline
$[p]_{8}$  & 3 & 1 &1  &  1  & 1 \\ \hline
$[p]_{9}$ & 4 & 1 & 1 & 1  & 1\\\hline
$[p]_{10}$ & 4 & 0 & 2 & 2  & 2\\\hline
$[p]_{11,\lambda}$ & 4 & 1 & 2 & 2  & 2\\\hline
\end{tabular}}
\end{table}
The result follows by Lemma \ref{ConjInvar} and Lemma \ref{p10conjugate}.
\end{proof}

\subsection{Case $L_{5}(\xi)$}

Recall that $L_{5}(\xi)$ for $\xi\in \mathbb{F}^{\times}$ has a basis $x, y, z, w$ with structure relations $[w,x]=x$ and $[w, y]=\xi y$. Clearly, $Z(L_{5}(\xi)) =\mspan\{z\}$. In addition, a linear map $\varphi:L_{5}(\xi)\to L_{5}(\xi)$ is a homomorphism of Lie algebras if and only if $\varphi(z)\in \mspan\{z\}$, $[\varphi(w),\varphi(x)]=\varphi(x)$ and $[\varphi(w),\varphi(y)]=\xi\varphi(y)$. It follows that $\varphi\in \mAut(L_{5}(\xi))$ if and only if it is given by
\begin{eqnarray}\label{isoformLe}
x\mapsto a_{1}x+a_{2} y, \quad y \mapsto b_{1}x+b_{2} y, \quad z\mapsto c_{3}z, \quad w\mapsto d_{1}x+d_{2}y+d_{3}z+d_{4}w,   
\end{eqnarray}
for some $a_{i}, b_{i}, c_{3},d_{i}\in \mathbb{F}$, which satisfy $(a_{1}b_{2}-a_{2}b_{1})c_{3}\neq 0$ and the following conditions:
\begin{eqnarray}\label{eqisoLe}
  a_{1}(d_{4}-1)=0, \quad a_{2}(\xi d_{4}-1)=0, \quad b_{1}(d_{4}-\xi)=0, \quad  b_{2}(d_{4}-1)=0.  
\end{eqnarray}
Note that if $\xi\neq-1$  then $d_{4}=1$, while if $\xi=-1$ then $d_{4}=\pm 1$.

\begin{prop}
There is no $p$-map on  $L_{5}(\xi)$ when $\xi \not\in \mathbb{F}_{p}^{\times}$. 
\end{prop}
\begin{proof}
Assume $\xi\not\in \mathbb{F}_{p}^{\times}$. Then $(\mad w)^{p} (x+y) =x+\xi^{p}y \not\in \mspan\{x+\xi y\}$. Thus $(\mad w)^{p}\neq \mad f$ for any $f\in L_{5}(\xi)$, because $[f, x+y] \in \mspan\{x+\xi y\}$. The result follows directly.   
\end{proof}

Assume $\xi\in \mathbb{F}_{p}^{\times}$. It is easy to verify that $(\mad x)^p= (\mad y)^p= (\mad z)^p= 0$ and $(\mad w)^p= \mad w$. Thus, by Lemma \ref{pmapconstruct}, the $p$-maps on $L_{5}(\xi)$ are the $p$-semilinear maps of the form
\begin{equation}
x^{[p]} =\alpha z, \quad y^{[p]} =\beta z, \quad z^{[p]} =\gamma z, \quad w^{[p]} =\delta z+w,
\end{equation}
where $\alpha, \beta, \gamma, \delta \in \mathbb{F}$. The next result classifies the $p$-maps on $L_{5}(\xi)$.

\begin{prop}\label{Lep}
Assume $\xi\in \mathbb{F}_{p}^{\times}$. The conjugacy classes of the $p$-maps on $L_{5}(\xi)$ are as follows: 
\begin{enumerate}
    \item $x^{[p]}=0,~y^{[p]}=0,~z^{[p]}=0,~w^{[p]}=w$;
    \item $x^{[p]}=z,~y^{[p]}=0,~z^{[p]}=0,~w^{[p]}=w$;
    \item $x^{[p]}=0,~y^{[p]}=z,~z^{[p]}=0,~w^{[p]}=w$, for $\xi\neq \pm1$; 
    \item $x^{[p]}=z,~ y^{[p]}=z,~z^{[p]}=0,~w^{[p]}=w$, for $\xi\neq 1$;
    \item $x^{[p]}=0,~y^{[p]}=0,~z^{[p]}=z,~w^{[p]}=w$;
    \item $x^{[p]}=z,~y^{[p]}=0,~z^{[p]}=z,~w^{[p]}=w$;
    \item $x^{[p]}=0,~y^{[p]}=z,~z^{[p]}=z,~w^{[p]}=w$, for $\xi\neq \pm1$;
    \item $x^{[p]}=z,~ y^{[p]}=z,~ z^{[p]}=z,~ w^{[p]}=w$, for $\xi\neq 1$.
\end{enumerate}
\end{prop}
\begin{proof}
Let $[p]$ be a $p$-map on $L_{5}(\xi)$. As explained above, we may assume $x^{[p]} =\alpha z,$ $y^{[p]} =\beta z,$ $z^{[p]} =\gamma z,$ and $ w^{[p]} =\delta z+w$ for some $\alpha, \beta, \gamma, \delta \in \mathbb{F}$. Let $[p]_{1}, \ldots, [p]_{8}$ be the $p$-maps listed in the proposition correspondingly. Ignore the restriction on $\xi$ in Parts (3), (4), (7) to (8) yet. Let $[p]':= \varphi\circ[p]\circ\varphi^{-1}$, where $\varphi\in \mAut(L_{5}(\xi))$ is the isomorphism given by
    \[
    x\mapsto x,\quad y\mapsto y,\quad z\mapsto z,\quad w\mapsto rz+w.
    \]
    Here, $r\in \mathbb{F}$ is chosen to satisfy $r+\delta-r^{p}\gamma=0$. It is easy to check that
    \[
    x^{[p]'}=\alpha z,\quad y^{[p]'}=\beta z,\quad z^{[p]'}=\gamma z,\quad w^{[p]'}=w.
    \]
According the  discussion  in Table \ref{arguLep}, 
\begin{table}[htbp]
\centering
\caption{Conjugate classes of $p$-maps of $L=L_{5}(\xi)$}\label{arguLep}
\setlength{\tabcolsep}{9mm}{
\renewcommand\arraystretch{2.5}
\begin{tabular}{c|c|c}
\hline
Parameters      & $\sigma=\sigma_{a,b,c,1}$ & $\sigma\circ[p]'\circ\sigma^{-1}$    \\ \hline
$\alpha=0$, $\beta=0$, $\gamma=0$  & $a=b=c=1$ & $[p]_{1}$ \\ \hline
$\alpha\neq0$, $\beta=0$, $\gamma=0$ & $a=\alpha^{\frac{1}{p}}$, $b=1$, $c=1$ & $[p]_{2}$    \\ \hline
$\alpha=0$, $\beta\neq0$, $\gamma=0$ & $a=1$, $b=\beta^{\frac{1}{p}}$, $c=1$ & $[p]_{3}$             \\ \hline
$\alpha\neq0$, $\beta\neq0$, $\gamma=0$ & $a=\alpha^{\frac{1}{p}}$, $b=\beta^{\frac{1}{p}}$, $c=1$ & $[p]_{4}$  \\\hline
$\alpha=0$, $\beta=0$, $\gamma\neq0$ & $a=1$, $b=1$, $c=\gamma^{\frac{1}{p-1}}$ & $[p]_{5}$  \\\hline
$\alpha\neq0$, $\beta=0$, $\gamma\neq0$ & $a=\left(\gamma^{\frac{1}{p-1}}\alpha\right)^{\frac{1}{p}}$, $b=1$, $c=\gamma^{\frac{1}{p-1}}$ & $[p]_{6}$    \\ \hline  
$\alpha=0$, $\beta\neq0$, $\gamma\neq0$ & $a=1$, $b=\left(\gamma^{\frac{1}{p-1}}\beta\right)^{\frac{1}{p}}$, $c=\gamma^{\frac{1}{p-1}}$ & $[p]_{7}$    \\  \hline  
$\alpha\neq0$, $\beta\neq0$, $\gamma\neq0$ & \makecell{$a=\left(\gamma^{\frac{1}{p-1}}\alpha\right)^{\frac{1}{p}}$, $b=\left(\gamma^{\frac{1}{p-1}}\beta\right)^{\frac{1}{p}}$ ,\\$c=\gamma^{\frac{1}{p-1}}$} & $[p]_{8}$  \\   \hline  
\end{tabular}}
\end{table}
$[p]'$  is conjugate to $[p]_{i}$ for some $i\in\{1,\ldots, 8\}$, and so is $[p]$. Furthermore, one has $\Delta_{\xi}\circ[p]_{3}\circ\Delta_{\xi}^{-1}=[p]_{2}$ and $\Delta_{\xi}\circ[p]_{7}\circ\Delta_{\xi}^{-1}=[p]_{6}$ for $\xi\in \{1, -1\}$,
where $\Delta_{\xi}\in \mAut(L_{5}(\xi))$ is the isomorphism given by 
\[
x\mapsto y,\quad y\mapsto x,\quad z\mapsto z,\quad w\mapsto \xi w;
\]
and $\psi\circ[p]_{4}\circ\psi^{-1}=[p]_{2}$ and $\psi\circ[p]_{8}\circ\psi^{-1}=[p]_{6}$ for $\xi=1$, where $\psi \in\mAut(L_{5}(1))$ is  given by 
\[
x\mapsto x,\quad y\mapsto x+y,\quad z\mapsto z,\quad w\mapsto w.
\]
As a summary, we show that $[p]$ is conjugate to one of the listed $p$-maps.

By a straightforward but tedious calculation on $\dim Z(L_{5}(\xi))^{[p]}$ and $\dim [L_{5}(\xi),L_{5}(\xi)]^{[p]}$, we have Table \ref{Lepdim}. 
\begin{table}[htbp]
\centering
\caption{2-tuple of invariant of $L=L_{5}(\xi)$}\label{Lepdim}
\setlength{\tabcolsep}{4.5mm}{
\renewcommand\arraystretch{1.6}
\begin{tabular}{ccccccccc}
\hline
$[p]$    & $[p]_{1}$ & $[p]_{2}$ & $[p]_{3}$ & $[p]_{4}$ & $[p]_{5}$ & $[p]_{6}$ & $[p]_{7}$ & $[p]_{8}$ \\ \hline
$\dim Z(L)^{[p]}$  & 0 & 0 & 0 & 0 & 1 & 1 & 1 & 1 \\ \hline
$\dim [L,L]^{[p]}$ & 0 & 1 & 1 & 1 & 0 & 1 & 1 & 1  \\ \hline
\end{tabular}} 
\end{table}
By Lemma \ref{ConjInvar}, it remains to show the following statements:
\begin{itemize}
    \item $[p]_{2}$ is not conjugate to $[p]_{4}$ for $\xi\neq 1$;
    \item $[p]_{6}$ is not conjugate to $[p]_{8}$ for $\xi\neq 1$;
    \item  $[p]_{3}$ is not conjugate to $[p]_{2}$ or $[p]_{4}$ for $\xi\neq \pm 1$;
    \item $[p]_{7}$ is not conjugate to $[p]_{6}$ or $[p]_{8}$ for $\xi\neq\pm 1$.
\end{itemize}
Suppose  there exists an isomorphism $\Phi\in\mAut(L_{5}(\xi))$ such that $\Phi\circ[p]_{2}\circ\Phi^{-1}=[p]_{4}$ (resp. $\Phi\circ[p]_{6}\circ\Phi^{-1}=[p]_{8}$) for some $\xi\neq 1$. Describe $\Phi(x), \Phi(y),\Phi(z),\Phi(w)$ as that in \eqref{isoformLe}.
Then 
\begin{align*}
    \Phi(x^{[p]_{2}})=\Phi(x)^{[p]_{4}} \quad (\text{resp. }\Phi(x^{[p]_{6}})=\Phi(x)^{[p]_{8}}) & \quad \Longrightarrow \quad c_{3}z=a_{1}^{p}z,\\
    \Phi(y^{[p]_{2}})=\Phi (y)^{[p]_{4}} \quad (\text{resp. }\Phi(y^{[p]_{6}})=\Phi(y)^{[p]_{8}}) & \quad\Longrightarrow \quad c_{3}z=b_{1}^{p}z.
\end{align*}
It follows that $a_{1}=b_{1}\neq0$. Put them into \eqref{eqisoLe}, one obtains $d_{4}=\xi=1$, which leads to a contradiction. Therefore, the first two desired statements hold. Finally, note that if $\xi\neq \pm 1$ then each isomorphism $\Psi\in\mAut(L_{5}(\xi))$ is necessarily of the form
\[
x\mapsto a_{1}x, \quad y \mapsto b_{2}y, \quad z\mapsto c_{3}z, \quad w\mapsto d_{1}x+d_{2}y+d_{3}z+w
\]
for some $a_{1},b_{2},c_{3}\in\mathbb{F}^{\times}$ and $d_{1},d_{2},d_{3}\in\mathbb{F}$. Thereof, if $\xi\neq \pm1$ then one may check that $\Psi\circ[p]_{3}\circ\Psi^{-1}(x)=0$ and $\Psi\circ[p]_{7}\circ\Psi^{-1}(x)=0$ for each $\Psi\in\mAut(L_{5}(\xi))$. Since $x^{[p]_{2}} =x^{[p]_{4}}=x^{[p]_{6}}=x^{[p]_{8}} =z$, it follows readily that the remaining two desired statements hold. 
\end{proof}

\begin{remark}\label{L5NONiso}
Let $[p]_{1}, \ldots, [p]_{8}$ be the $p$-maps on $L_{5}(\xi)$ listed in Proposition \ref{Lep}. Recall that $\Xi_{p}:=\{\xi_{p}^{r} ~ |~ 0\leq r\leq \frac{p-1}{2}\}$, where $\xi_{p} \in \mathbb{F}_{p}^{\times}$ is a chosen generator of $\mathbb{F}_{p}^{\times}$. Since $$\Xi_{p}\cup \Xi_{p}^{-1} = \mathbb{F}_{p}^{\times} \quad \text{and} \quad \Xi_{p}\cap \Xi_{p}^{-1} = \{\pm 1\},$$  Proposition \ref{Liewith3ideal} and Proposition \ref{Lep} tell us that the family of restricted Lie algebras $(L_{5}(\xi), [p]_{i})$, where $\xi\in \Xi_{p}$ and $i\in \{1,\ldots, 8\}$, are mutually non-isomorphic.

Furthermore, for any $\xi\in\mathbb{F}^{\times}_{p}\backslash \{\pm 1\}$,  there exists a unique permutation $\tau\in S_{8}$ such that  $$(L_{5}(\xi), [p]_{\tau(i)}) \cong (L_{5}(\xi^{-1}), [p]_{i})$$ as restricted Lie algebras. It turns out that $\tau =(23) (67)$. 
To see this, let $\varphi: L_{5}(\xi) \to L_{5}(\xi^{-1})$ be an arbitrary isomorphism of Lie algebras. One candidate is the one given by $$x\mapsto y, \quad y\mapsto x, \quad z\mapsto z, \quad  w\mapsto \xi w.$$
Then $[p]_{\tau(i)}$  and $\varphi^{-1} \circ [p]_{i}\circ \varphi$ are conjugate as $p$-maps on $L_{5}(\xi)$ for each $i=1,\ldots,8$. One may obtain the value of $\tau(i)$ by a tedious but straightforward calculation.
\end{remark}

\subsection{Case $L_{6}(\xi,\eta)$}

Recall that the Lie algebra $L_{6}(\xi, \eta)$ for $\xi, \eta\in \mathbb{F}^{\times}$ has a basis $x, y, z, w$ with structure relations $[w,x]=x$, $[w, y]=\xi y$ and $[w, z]=\eta z$. Clearly, the center is trivial.

\begin{prop}
There is no $p$-map on $L_{6}(\xi,\eta)$ when $(\xi,\eta)\not\in \mathbb{F}_{p}^{\times} \times \mathbb{F}_{p}^{\times}$. But there is exactly one $p$-map on $L_{6}(\xi,\eta)$ for any $(\xi,\eta)\in \mathbb{F}_{p}^{\times} \times \mathbb{F}_{p}^{\times}$, which is
\begin{equation*}
     x^{[p]}=0, \quad y^{[p]}=0, \quad z^{[p]}=0, \quad w^{[p]}=w.
\end{equation*}
\end{prop}
\begin{proof} 
If  $(\xi,\eta)\not\in \mathbb{F}_{p}^{\times} \times \mathbb{F}_{p}^{\times}$ then $(\mad w)^{p} (x+y+z) =x+\xi^{p}y +\eta^{p}z \not\in \mspan\{x+\xi y+\eta z\}$ and hence  $(\mad w)^{p}\neq \mad f$ for any $f\in L_{6}(\xi,\eta)$, because $[f, x+y+z] \in \mspan\{x+\xi y+\eta z\}$. The first statement follows directly. Otherwise, one has $(\mad x)^p = (\mad y)^p =(\mad z)^p =0$ and $(\mad w)^p = \mad w$.  Since the center is trivial, the second statement follows by Lemma \ref{pmapconstruct}.
\end{proof}

\begin{remark}\label{L6NONiso}
Let $\xi, \eta, \xi', \eta' \in \mathbb{F}_{p}^{\times}$. Let $[p]$ and $[p]'$ be the unique $p$-map on $L_{6}(\xi, \eta)$ and $L_{6}(\xi', \eta')$ respectively. Clearly, $(L_{6}(\xi, \eta), [p]) \cong (L_{6}(\xi', \eta'), [p]')$ as restricted Lie algebras is equivalent to that  $L_{6}(\xi, \eta) \cong  L_{6}(\xi', \eta')$ as Lie algebras. So  by  Proposition \ref{Liewith3ideal},  $(L_{6}(\xi, \eta), [p]) \cong (L_{6}(\xi', \eta'), [p]')$ if and only if   $(\xi, \eta) = \tau \cdot (\xi', \eta')$ for some permutation $\tau \in S_{3}$.
\end{remark}

\subsection{Case $L_{7},L_{8}(\xi),L_{9}$}

Recall that the structure relations of $L_{7}$ are $[w,x]=y$ and $[w,z]=z$; and those of $L_{8}(\xi)$ (resp. $L_{9}$) are $[w,x]=x+y$, $[w,y]=y$ and $[w,z]=\xi z$ (resp. $[w,z]=x+z$).

\begin{prop}\label{norestric}
There is no $p$-map on $L_{7}$, $L_{8}(\xi)$ and $L_{9}$.
\end{prop}
\begin{proof}
It is easy to check that $(\mad w)^{p}(x+z)=z$ and $[f, x+z] \in \mspan\{ y+z\}$ in $L_{7}$. Also, one readily obtains that  $(\mad w)^{p}(x) =x$ and $[f,x] \in \mspan\{x+y\}$ in $L_{8}(\xi)$ and $L_{9}$. Thus $(\mad w)^{p}\neq \mad f$ for any $f$ in $L_{7}$, $L_{8}(\xi)$ or $L_{9}$. The result follows readily. 
\end{proof}

\section{The $p$-maps on solvable Lie algebras of dimension $4$ without $3$-dimensional abelian ideals}
\label{Classify-restrictedLie-2}

This section is devoted to classify  the $p$-maps on the Lie algebras presented in Proposition \ref{Liewithout3ideal}. Note that these Lie algebras  form a complete list of representatives of solvable Lie algebras of dimension $4$ without $3$-dimensional abelian ideals.

\subsection{Case $N_{1}$}

Recall that the Lie algebra $N_{1}$ has a basis $x, y, z, w$ with structure relations $[y,x]=x$ and $[w,z]=z$. Clearly, the center of $N_{1}$ is trivial.

\begin{prop}\label{Lgp}
 There is exactly one $p$-map on $N_{1}$, which is
\begin{equation*}
    x^{[p]}=0, \quad  y^{[p]}=y, \quad z^{[p]}=0, \quad w^{[p]}=w.
\end{equation*}
\end{prop}
\begin{proof}
It is easy to verify that $(\mad x)^p =(\mad z)^p =0$, $(\mad y)^{p} =\mad y$ and $(\mad w)^p = \mad w$.  Since the center of  $N_{1}$ is trivial, the result follows by Lemma \ref{pmapconstruct}.
\end{proof}

\subsection{Case $N_{2}$}

Recall that the Lie algebra $N_{2}$ has a basis $x, y, z, w$ with structure relations  $[z,x]=y$, $[w,x]=x$, $[w,y]=2y$ and $[w,z]=z$. Clearly, if $p=2$ then  $Z(N_{2}) =\mspan\{ y\}$, while if $p\geq 3$ then $Z(N_{2})=0$. In addition, a linear map $\varphi:N_{2}\mapsto N_{2}$ is a homomorphism of Lie algebras if and only if $\varphi(y)\in Z(N_{2}) =\mspan\{ y\}$,  $[\varphi(z),\varphi(x)]=\varphi(y)$, $[\varphi(w),\varphi(x)]=\varphi(x)$, $[\varphi(w),\varphi(y)]=2\varphi(y)$ and $[\varphi(w),\varphi(z)]=\varphi(z)$. It follows that $\varphi\in\mAut(N_{2})$ if and only if it is given by

\begin{equation}
\begin{aligned}
x&\mapsto a_{1}x+(a_{3}d_{1}- a_{1}d_{3})y+a_{3}z,& y&\mapsto (a_{1}c_{3}-a_{3}c_{1}) y,\\ 
z&\mapsto c_{1}x+(c_{3}d_{1}-c_{1}d_{3})y+c_{3}z,& w&\mapsto d_{1}x+d_{2}y+d_{3}z+w
\end{aligned}
\end{equation}
for some $a_{i},c_{i},d_{i}\in \mathbb{F}$, which satisfy $a_{1}c_{3}-a_{3}c_{1}\neq0$.

Assume $p=2$. It is easy to verify that $(\mad x)^2 =(\mad z)^2 =(\mad y)^{2} =0$ and $(\mad w)^2 = \mad w$. Thus by Lemma \ref{pmapconstruct}, the $2$-maps on $N_{2}$ are the $2$-semilinear maps of the form:
\begin{eqnarray}
x^{[2]}=\alpha y,\quad y^{[2]}=\beta y,\quad z^{[2]}=\gamma y,\quad w^{[2]}=\delta y+w,
\end{eqnarray}
where $\alpha, \beta, \gamma, \delta\in \mathbb{F}$. The next result classifies the $2$-maps on $N_{2}$.

\begin{prop}\label{Lhp}
Let $p=2$. The conjugacy classes of the $2$-maps on $N_{2}$ are as follows:
\begin{enumerate}
    \item $x^{[2]}=0,~ y^{[2]}=0,~ z^{[2]}=0,~ w^{[2]}=w$.
    \item $x^{[2]}=y,~ y^{[2]}=0,~ z^{[2]}=0,~ w^{[2]}=w$.
    \item $x^{[2]}=0,~ y^{[2]}=y,~ z^{[2]}=0,~ w^{[2]}=w$.
\end{enumerate}

\end{prop}
\begin{proof}
Let $[2]$ be an arbitrary $2$-map on $N_{2}$. As explained above, we may assume $x^{[2]}=\alpha y,~ y^{[2]}=\beta y,~ z^{[2]}=\gamma y,~ w^{[2]}=\delta y+w$ for some $\alpha, \beta, \gamma, \delta\in \mathbb{F}$. Let $[2]_{1}, [2]_{2}, [2]_{3}$ be the $2$-maps listed in the proposition correspondingly. We split the discussion as follows.
    \begin{enumerate}
    \item Assume $\beta=0$ and $\alpha=\gamma=0$. Then $\psi\circ [p]\circ\psi^{-1}=[2]_{1}$, where  $\psi\in \mAut(N_{2})$ is the isomorphism given by
    \[
    x\mapsto x,\quad y\mapsto  y,\quad z\mapsto z,\quad w\mapsto \delta y+w.
    \]
    \item  Assume $\beta=0$, $\alpha=0$ and $\gamma\neq0$. Then $\chi \circ [2] \circ \chi^{-1} =[2]_{2}$, where $\chi \in \mAut(N_{2})$ is the isomorphism given by
        \[
        x\mapsto \gamma^{\frac{1}{2}} z, \quad y\mapsto y, \quad z\mapsto \gamma^{-\frac{1}{2}} x, \quad  w\mapsto \delta y+ w.
        \]
    \item Assume $\beta=0$ and $\alpha\neq0$. Then $\tau \circ [2]\circ \tau^{-1}=[2]_{2}$, where $\tau\in \mAut(N_{2})$ is the isomorphism given by 
        \[
        x\mapsto \alpha^{\frac{1}{2}}x, \quad  y\mapsto y, \quad z\mapsto \gamma^{\frac{1}{2}}x+ \alpha^{-\frac{1}{2}}z, \quad w\mapsto \delta y+ w.
        \]
    \item Assume  $\beta\neq0$. Then $\varphi\circ[2]\circ\varphi^{-1}=[2]_{3}$, where $\varphi\in\mAut(N_{2})$ is the isomorphism given by 
    \[
    x\mapsto \beta^{\frac{1}{2}} x+(\alpha\beta)^{\frac{1}{2}}y,\quad y\mapsto \beta y,\quad z\mapsto (\gamma\beta )^{\frac{1}{2}}y+\beta^{\frac{1}{2}}z,\quad w\mapsto \gamma^{\frac{1}{2}}x+d_{2}y+\alpha^{\frac{1}{2}}z+w.
    \]
    Here, $d_{2}\in \mathbb{F}$ is any scalar that satisfies $d_{2}^{2}-d_{2}=\delta\beta$.
    \end{enumerate}
    
It remains to show $[2]_{1}, [2]_{2}, [2]_{3}$  are not conjugate to each other.
Note that $Z(N_{2})^{[2]_{1}}=Z(N_{2})^{[2]_{2}}=0$ but $Z(N_{2})^{[2]_{3}}\neq0$, so $[2]_{3}$ is not conjugate to $[2]_{1}$ and $[2]_{2}$.
In addition, $[N_{2},N_{2}]^{[2]_{1}}=0$ but $[N_{2},N_{2}]^{[2]_{2}}\neq 0$, it follows that $[2]_{1}$is not conjugate to $[2]_{2}$.
\end{proof}

\begin{prop}
 Let $p\geq 3$. There is exactly one $p$-map on $N_{2}$, which is 
 \[
x^{[p]}=0,\quad y^{[p]}=0,\quad z^{[p]}=0,\quad w^{[p]}=w.
 \]
\end{prop}
\begin{proof}
    It is straightforward to verify that $(\mad x)^p =(\mad z)^p=(\mad y)^{p} =0$ and $(\mad w)^p = \mad w$.  Since the center of  $N_{2}$ is trivial when $p\geq 3$, the result follows by Lemma \ref{pmapconstruct}.
\end{proof}

\subsection{Case $N_{3}(\xi)$}

Recall that the Lie algebra $N_{3}(\xi)$, where $\xi\in \mathbb{F}$, has a basis $x, y, z, w$ with structure relations $[z,x]=y$, $[w,x]=x+\xi z$, $[w,y]=y$ and $[w,z]=x$. Clearly, the center is trivial.

\begin{prop}
Let  $p=2$. There is no  $2$-map on $N_{3}(\xi)$ when $\xi\neq 0$. But there is exactly one $2$-map on $N_{3}(0)$, which is 
\[
x^{[2]}=0,\quad y^{[2]}=0,\quad z^{[2]}=y, \quad w^{[2]}=w.
\]
\end{prop}

\begin{proof}
If $\xi\neq 0$ then $(\mad w)^{2}(z) =x+\xi z\not\in \mspan\{x,y\}$ and hence $(\mad w)^{2}\neq \mad f$ for any $f\in N_{3}(\xi)$, because $[f, z]\in \mspan\{x,y\}$. The first statement follows directly. In addition, it is easy to verify that  $(\mad x)^{2} =(\mad y)^{2} = 0$, $(\mad z)^{2} = \mad y$ and $(\mad w)^{2} =\mad w$ when $\xi =0$. Since the center is trivial, the second statement follows by Lemma \ref{pmapconstruct}.
\end{proof}

Recall that $Q_{p}=\{ \alpha\in \mathbb{F}_{p}^{\times} ~|~ \alpha^{\frac{p-1}{2}}=1 \}$. It is the set of quadratic residues modulo $p$. 

\begin{prop}
Let $p\geq 3$. There  is no $p$-map on $N_{3}(\xi)$ when $\xi\not\in Q_{p}-\frac{1}{4}$. But there is exactly one $p$-map on $N_{3}(\xi)$ for any $\xi\in Q_{p}-\frac{1}{4}$, which is 
\[
x^{[p]}=0,\quad y^{[p]}=0,\quad z^{[p]}=0, \quad w^{[p]}=w.
\]
\end{prop}

\begin{proof}
It  is easy to verify that $(\mad x)^{p}= (\mad y)^{p}= (\mad z)^{p}= 0$, $(\mad w)^{p}(y) =y$ and $(\mad w)^{p}(w) =0$. By using the Jordan normal form of square matrices, one has
\begin{eqnarray}\label{adwLi}
\big ( (\mad w)^{p}(x), (\mad w)^{p} (z) \big ) = (x, z)  \begin{pmatrix}[1.3]
        1& 1\\ \xi &0
    \end{pmatrix}^{p}
    =  (x, z)  \begin{pmatrix}[1.3]
        \frac{1}{2} + \frac{1}{2}\eta& \eta \\ \xi \eta & \frac{1}{2} - \frac{1}{2}\eta
    \end{pmatrix},
\end{eqnarray}
where $\eta:= (1+4\xi)^{\frac{p-1}{2}}$. Note that $\eta=1$ if and only if $\xi\in Q_{p}-\frac{1}{4}$. If $\eta\neq 1$ then  $(\mad w)^{p}(z) = \eta x+\frac{1}{2}(1-\eta) z \not \in \mspan\{x, y\}$ by \eqref{adwLi}, and hence $(\mad w)^{p}\neq \mad f$ for any $f\in N_{3}(\xi)$, because $[f, z]\in \mspan\{x, y\}$. The first statement follows directly. Otherwise, if  $\eta=1$ then  $(\mad w)^{p} =\mad w$ by \eqref{adwLi}. Since the center is trivial, the second statement follows by Lemma \ref{pmapconstruct}. 
\end{proof}

\subsection{Case $N_{4}$}

Recall that the Lie algebra $N_{4}$ has a basis $x, y, z ,w$ with structure relations $[z,x]=y$, $[w,x]=z$ and $[w,z]=x$. Clearly, $Z(N_{4}) =\mspan\{y\}$.  In addition, a linear map $\varphi:N_{4}\mapsto N_{4}$ is a homomorphism of Lie algebras if and only if $\varphi(y)\in \mspan\{y\}$, $[\varphi(z),\varphi(x)]=\varphi(y)$, $[\varphi(w),\varphi(x)]=\varphi(z)$ and $[\varphi(w),\varphi(z)]=\varphi(x)$. It follows that $\varphi\in\mAut(N_{4})$ if and only if it is given by
\begin{equation}\label{isoformLj}
    \begin{aligned}
    x&\mapsto a_{1}x+ (a_{1}d_{3}-c_{1}d_{1}d_{4})y+d_{4}c_{1}z, \quad  \quad \quad&
    y&\mapsto (a_{1}^2-c_{1}^2)d_{4} y,\\
    z&\mapsto c_{1}x+(c_{1}d_{3}-a_{1}d_{1}d_{4})y+d_{4}a_{1}z, \quad \quad \quad & 
    w&\mapsto d_{1}x+d_{2}y+d_{3}z+d_{4}w
    \end{aligned}
    \end{equation}
for some $a_{1},c_{1},d_{1}, d_{2}, d_{3}, d_{4}\in \mathbb{F}$, which satisfy $d_{4}^2=1$ and $a_{1}^2-c_{1}^2\neq0$. 

\begin{prop}
Let $p=2$. There is no $2$-map on $N_{4}$.
\end{prop}
\begin{proof}
Clearly, $(\mad w)^{2}(x)=x$. Thus $(\mad w)^{2}\neq \mad f$ for any $f\in N_{4}$, because $[f,x]\in \mspan\{y,z\}$. The result follows directly. 
\end{proof}

Assume $p\geq 3$. It is easy to verify that $(\mad x)^{p} = (\mad y)^{p} = (\mad z)^{p}=0$ and $ (\mad w)^{p} =\mad w$. Thus by Lemma \ref{pmapconstruct}, the $p$-maps on $N_{4}$ are the $p$-semilinear maps of the form:
\[
x^{[p]}=\alpha y, \quad  y^{[p]}=\beta y, \quad z^{[p]}=\gamma y, \quad w^{[p]}=\delta y +w,
\]
where $\alpha, \beta, \gamma, \delta \in \mathbb{F}$. The following result classifies the $p$-maps on $N_{4}$ for $p\geq 3$.

\begin{prop}\label{Ljp}
Let $p\geq3$. The conjugacy classes of the $p$-maps on $N_{4}$ are as follows:
\begin{enumerate}
    \item $x^{[p]}=0,~y^{[p]}=0,~z^{[p]}=0,~w^{[p]}=w$;
    \item $x^{[p]}=y,~y^{[p]}=0,~z^{[p]}=0,~w^{[p]}=w$;
    \item $x^{[p]}=0,~y^{[p]}=y,~z^{[p]}=0,~w^{[p]}=w$;
    \item $x^{[p]}=y,~y^{[p]}=0,~z^{[p]}=y,~w^{[p]}=w$.
\end{enumerate}
\end{prop}

To justify the proposition, we need the following lemma.

\begin{lemma}\label{pstructureofLj}
    Let $p\geq3$ and $[p]$ be a $p$-semilinear map on $N_{4}$. Then for any $a,b,c,d\in\mathbb{F}$,
    \[
    (ax+by+cz+dw)^{[p]}=a^{p}x^{[p]}+b^{p}y^{[p]}+c^{p}z^{[p]}+d^{p}w^{[p]}+d^{p-1}(ax+cz)-\frac{1}{2}d^{p-2}(a^{2}-c^{2})y.
    \]
\end{lemma}
\begin{proof}
    Since $y$ belongs to the center of $N_{4}$ it suffices to show
    \[
    (ax+cz+dw)^{[p]}=a^{p}x^{[p]}+c^{p}z^{[p]}+d^{p}w^{[p]}+d^{p-1}(ax+cz)-\frac{1}{2}d^{p-2}(a^{2}-c^{2})y.
    \]
    First, we show that
    \[
    (\mad (dwT+ax+cz))^{k}(dw)=-d^{k}(ax+cz)T^{k-1}+d^{k-1}(a^{2}-c^{2})yT^{k-2}    \]
    for even integer $k\geq0$. We prove it by induction on $k$. When $k=2$, one has
    \begin{align*}
    [dwT+ax+cz,[dwT+ax+cz,dw]]&=[dwT+ax+cz,-d(cx+az)]\\
    &=-d^{2}(ax+cz)T+d(a^{2}-c^{2}).
    \end{align*}
    Then by the induction hypothesis, one has
    \begin{align*}
        (\mad (dwT+ax+cz))^{k}(dw)
        &= [dwT+ax+cz,[dwT+ax+cz,(\mad (dwT+ax+cz))^{k-2}(dw)]]\\
        &=[dwT+ax+cz,-d^{k-1}(cx+az)T^{k-2}]\\
        &=-d^{k}(ax+cz)T^{k-1}+d^{k-1}(a^{2}-c^{2})yT^{k-2}.
    \end{align*}
    Hence, by the equation \eqref{adT-1}, one has
    \begin{align*}
        (ax+cz+dw)^{[p]}
        &=a^{p}x^{[p]}+c^{p}z^{[p]}+d^{p}w^{[p]}-\frac{d^{p-1}}{p-1}(ax+cz)+\frac{d^{p-2}}{p-2}(a^{2}-c^{2})y\\
        &=a^{p}x^{[p]}+c^{p}z^{[p]}+d^{p}w^{[p]}+d^{p-1}(ax+cz)-\frac{d^{p-2}}{2}(a^{2}-c^{2})y
    \end{align*}
    as desired.
\end{proof}

\begin{proof}[Proof of Proposition \ref{Ljp}]
Let $[p]$ be an arbitrary $p$-map on $N_{4}$. As explained above, we may assume $x^{[p]}=\alpha y$, $y^{[p]}=\beta y$, $z^{[p]}=\gamma y$ and $w^{[p]}=\delta y +w$ for some $\alpha, \beta, \gamma, \delta \in \mathbb{F}$. Let $[p]'= \varphi\circ[p]\circ\varphi^{-1}$, where $\varphi\in \mAut(N_{4})$ is the isomorphism given by
\[x\mapsto x,\quad y\mapsto y,\quad z\mapsto z,\quad w\mapsto uy +w.\]
Here, $u\in\mathbb{F}$ is any scalar that satisfies $\delta-u^{p}\beta+u=0$. Readily, one has 
\[
x^{[p]'}=\alpha y,\quad y^{[p]'}=\beta y,\quad z^{[p]'}=\gamma y,\quad w^{[p]'}=w.
\]
Let $[p]_{1}, [p]_{2}, [p]_{3}, [p]_{4}$ be the $p$-maps listed in the proposition correspondingly. We proceed to show that $[p]'$ is conjugate to one of them as follows.
\begin{enumerate}
    \item Assume $\gamma=0$, $\beta=0$ and $\alpha=0$. Then $[p]'=[p]_{1}$.
    \item Assume $\gamma=0$,  $\beta=0$ and $\alpha\neq 0$. Then $\sigma_{a,b,c,1}\circ[p]'\circ\sigma_{a,b,c,1}^{-1}=[p]_{2}$, where $a=c=\alpha^{\frac{1}{p-2}}$ and $b=\alpha^{\frac{2}{p-2}}$.
    \item  Assume  $\gamma=0$ and $\beta\neq 0$. Let $\varphi\in \mAut(N_{4})$ be the isomorphism given by
\begin{align*}
    x&\mapsto\beta^{\frac{1}{2(p-1)}}x+\beta^{\frac{1}{p(p-1)}}\alpha^{\frac{1}{p}}y,&
    y&\mapsto\beta^{\frac{1}{p-1}}y,\\
    z&\mapsto\beta^{\frac{1}{2(p-1)}}z,& 
    w&\mapsto -\beta^{\frac{1}{2(p-1)}}\left(\frac{\alpha}{\beta}\right)^{\frac{1}{p}}x-vy+w,
\end{align*}
where $v\in \mathbb{F}$ is any scalar that satisfies $v^{p}-v=\frac{1}{2}(\frac{\alpha}{\beta})^{\frac{2}{p}}\beta^{\frac{1}{p-1}}$. By Lemma \ref{pstructureofLj}, it is easy to verify $\varphi\circ[p]\circ\varphi^{-1}=[p]_{3}$.

  \item  Assume $\gamma\neq 0$ and $\alpha=0$. Let $[p]'':=\Delta\circ[p]'\circ\Delta^{-1}$, where $\Delta\in \mAut(N_{4})$ is given by
    \[
    x\mapsto -z,\quad y\mapsto -y,\quad z\mapsto -x,\quad w\mapsto w.
    \]
    Readily, one has
    \[
    x^{[p]''}=\gamma y,\quad y^{[p]''}=\beta y,\quad z^{[p]''}=\alpha y=0,\quad w^{[p]''}=w.
    \]
   Hence, $[p]''$  falls into the situations (2) or (3).
   
   \item Assume $\gamma,\alpha\neq 0$ and $\alpha\neq\pm\gamma$. Let $\psi\in \mAut(N_{4})$ be the isomorphism given by
\[
x\mapsto \frac{1}{1-t^{2}}x-\frac{t}{1-t^{2}}z,\quad y\mapsto \frac{1}{1-t^{2}}y,\quad z\mapsto -\frac{t}{1-t^{2}}x+\frac{1}{1-t^{2}}z,\quad w\mapsto w,
\]
where $t=(\frac{\gamma}{\alpha})^{\frac{1}{p}}$. Then $[p]'':=\psi\circ [p]'\circ \psi^{-1}$ satisfies that
\[
x^{[p]''}=\frac{\alpha^{2}+\gamma^{2}}{\alpha(1-t^{2})}y,\quad y^{[p]''}=(1-t^{2})^{p-1}\beta y,\quad z^{[p]''}=0,\quad w^{[p]''}=w.
\]
Thus $[p]''$ falls into the situations (1), (2) or (3).

\item Assume $\alpha=\pm\gamma\neq 0$ and $\beta=0$. Then $\sigma_{a,a^{2},\pm a,\pm 1}\circ[p]'\circ\sigma_{a,a^{2},\pm a,\pm 1}^{-1}=[p]_{4}$, where $a=\alpha^{\frac{1}{p}}$.

\item Assume $\alpha=\pm\gamma\neq 0$ and $\beta\neq0$. Let $\psi\in\mAut(N_{4})$ be the isomorphism given by
\begin{align*}
    x&\mapsto \frac{r\pm s}{2}x + ry +\frac{r\mp s}{2}z,& y&\mapsto rsy\\
    z&\mapsto \frac{\pm r-s}{2}x\pm ry+\frac{\pm r+s}{2}z & w&\mapsto \mp x\pm z\pm w,
\end{align*}
where $s=(\frac{\beta}{\alpha})^{\frac{1}{p}}$ and $r=(\alpha s)^{\frac{1}{p-1}}$. By Lemma \ref{pstructureofLj}, it is easy to verify $\psi\circ[p]\circ\psi^{-1}=[p]_{3}$.
\end{enumerate}

By a straightforward but tedious calculation on $\dim Z(N_{4})^{[p]}$ and $\dim [N_{4},N_{4}]^{[p]}$, we have Table \ref{Ljpdim}
\begin{table}[htbp]
\centering
\caption{2-tuple of invariant of $L=N_{4}$}\label{Ljpdim}
\setlength{\tabcolsep}{4.0mm}{
\renewcommand\arraystretch{1.4}
\begin{tabular}{ccccc}
\hline
$[p]$      & $[p]_{1}$ & $[p]_{2}$ & $[p]_{3}$ & $[p]_{4}$\\ \hline
$\dim Z(L)^{[p]}$  & 0 & 0 & 1 & 0   \\ \hline
$\dim [L,L]^{[p]}$ & 0 & 1 & 1 & 1  \\ \hline
\end{tabular}} 
\end{table}
So by Lemma \ref{ConjInvar}, it remains to show that $[p]_{2}$ and $[p]_{4}$ are not conjugate.
Suppose there exists an isomorphism $\Phi\in\mAut(N_{4})$ such that $\Phi\circ[p]_{4}\circ\Phi^{-1}=[p]_{2}$. Describe $\Phi(x), \Phi(y),\Phi(z),\Phi(w)$ as that in \eqref{isoformLj}. Then one has the following implications:
\begin{align*}
    \Phi(x^{[p]_{4}})&=\Phi(x)^{[p_{2}]} \quad \Longrightarrow \quad b_{2}=a_{1}^{p},\\
    \Phi(z^{[p]_{4}})&=\Phi(z)^{[p_{2}]}\quad\Longrightarrow \quad b_{2}=c_{1}^{p},
\end{align*}
It follows that $a_{1}=c_{1}$, which leads to a contradiction.
\end{proof}

\subsection{Case $N_{5}$}

Recall that the Lie structure  $N_{5}$ has a basis $x,y,z,w$ with structure relations $[z,x]=y$,  $[z,y]=x$, $[w,x]=x$ and $[w,z]=z$. One need to assume $p=2$ for $N_{5}$ (see Remark \ref{N5p3}).

\begin{prop}
Let $p=2$. There is no $2$-map on $N_{5}$.
\end{prop}
\begin{proof}
Clearly, $(\mad z)^{2}(y) =y$. Thus $(\mad z)^{2}\neq \mad f$ for any $f\in N_{5}$, because $[f,y]\in \mspan\{x\}$. The result follows directly.
\end{proof}

\section{The $p$-maps on non-solvable Lie algebras of dimension $4$}

\label{Classify-restrictedLie-3}

This section focus on the classification of the $p$-maps on the Lie algebras presented in Proposition \ref{nonsolvLie}, which form a complete list of representatives  of non-solvable Lie algebras of dimension $4$.

\subsection{Even characteristic}

 Assume $p=2$. By  Proposition \ref{nonsolvLie}, there are exactly two non-solvable Lie algebras of dimension $4$ up to isomorphism. They are $W_{1}$ and $W_{2}$ in our notation, and both have  $x, y, z, w$ as a basis. The structure relations of $W_{1}$ are $[y,x]=z$, $[z,x]=x$ and $[z,y]=y$; while those of $W_{2}$ are $[y,x]=z$, $[z,x]=x$, $[z,y]=y$ and $[w,x]=y$.

 \begin{prop}
Let $p=2$. There is no $2$-map on  $W_{1}$ and $W_{2}$. 
 \end{prop}
\begin{proof}
Let $L$ be $W_{1}$ or $W_{2}$. Clearly, $(\mad x)^{2}(y)=x$. Thus $(\mad x)^{2}\neq \mad f$ for any $f\in L$, because $[f, y]\in \mspan\{ y, z\}$. The result follows directly. 
\end{proof}

\subsection{Odd characteristic}

Assume $p\geq 3$. By  Proposition \ref{nonsolvLie},  there is only one non-solvable Lie algebra of dimension $4$ up to isomorphism. It is $\mathfrak{gl}_2$ in our notation, and has a basis $x, y, z, w$ with structure relations 
\[
[y,x] =-z, \quad [z,x] =2x, \quad [z,y] =-2y.
\]
Clearly, the center is $Z(\mathfrak{gl}_2)=\mspan\{w\}$. In addition, a linear map $\varphi:\mathfrak{gl}_2\to \mathfrak{gl}_2$ is a homomorphism of Lie algebras if and only if  $\varphi(w) \in \mspan\{w\}$, $[\varphi(y),\varphi(x)] = -\varphi(z)$, $[\varphi(z),\varphi(x)] =2\varphi(x)$ and $[\varphi(z),\varphi(y)] =2 \varphi(y)$. It follows that $\varphi \in \mAut(\mathfrak{gl}_2)$ if and only if it is given by 
\begin{eqnarray}\label{AutoLk}
  x\mapsto a_{1}x+a_{2}y+a_{3}z, \quad y\mapsto b_{1}x+b_{2}y+b_{3}z, \quad z\mapsto c_{1}x+c_{2}y+c_{3}z, \quad w\mapsto d_{4}w    
\end{eqnarray}
for some  $a_{i},b_{i},c_{i},d_{4}\in\mathbb{F}$, which satisfy the following conditions:
\begin{equation}\label{eqisogl2}
    \begin{aligned}
         a_{1}c_{3}-a_{3}c_{1}-a_{1}&=0,\quad& 
         a_{3}c_{2}-a_{2}c_{3}-a_{2}&=0,\quad& 
         a_{2}c_{1}-a_{1}c_{2}-2a_{3}&=0,\\
         b_{3}c_{1}-b_{1}c_{3}-b_{1}&=0,\quad&
         b_{2}c_{3}-b_{3}c_{2}-b_{2}&=0,\quad& 
         b_{1}c_{2}-b_{2}c_{1}-2b_{3}&=0,\\
         2a_{3}b_{1}-2a_{1}b_{3}-c_{1}&=0,\quad& 
         2a_{2}b_{3}-2a_{3}b_{2}-c_{2}&=0,\quad& 
         a_{1}b_{2}-a_{2}b_{1}-c_{3}&=0,
    \end{aligned}
    \end{equation}
as well as $(a_{1}b_{2}c_{3}+a_{2}b_{3}c_{1}+a_{3}b_{1}c_{2}-a_{1}b_{3}c_{2}-a_{2}b_{1}c_{3}-a_{3}b_{2}c_{1}) d_{4} \neq0.$

It is straightforward to verify that $(\mad x)^{p} = (\mad y)^{p} = (\mad w)^{p} =0$ and $(\mad z)^p =\mad z$ for $\mathfrak{gl}_2$. Thus by Lemma \ref{pmapconstruct}, the $p$-maps on $\mathfrak{gl}_2$ are the $p$-semilinear maps of the form:
\begin{eqnarray}
x^{[p]}=\alpha w, \quad y^{[p]}=\beta w, \quad z^{[p]}=z+\gamma w, \quad w^{[p]}=\delta w,  
\end{eqnarray}
where $\alpha, \beta, \gamma, \delta\in \mathbb{F}$. The next result classifies the $p$-maps on $\mathfrak{gl}_2$.

\begin{prop}\label{Lkp}
Let $p\geq 3$. The conjugacy classes of the $p$-maps on $\mathfrak{gl}_2$ are as follows:
    \begin{enumerate}
        \item $x^{[p]}=0,~y^{[p]}=0,~z^{[p]}=z,~w^{[p]}=0$;
        \item $x^{[p]}=w,~y^{[p]}=0,~z^{[p]}=z,~w^{[p]}=0$;
        \item $x^{[p]}=0,~y^{[p]}=0,~z^{[p]}=z,~w^{[p]}=w$;
        \item $x^{[p]}=w,~y^{[p]}=0,~z^{[p]}=z,~w^{[p]}=w$;
        \item $x^{[p]}=0, ~ y^{[p]}=0,~z^{[p]}=z+w,~w^{[p]}=\lambda w$, $\lambda\in\mathbb{F}$.
        \end{enumerate}
\end{prop}

The proof of this proposition will be addressed at the end of this subsection. In the remainder of this subsection, we tacitly assume $p\geq 3$ and denote by $[p]_1$, $[p]_2$, $[p]_3$, $[p]_4$, $[p]_{5,\lambda}$, where $\lambda\in \mathbb{F}$,  the $p$-maps given in Proposition \ref{Lkp} correspondingly. 

\begin{remark}
Note that $\mathfrak{gl}_2 = \mathfrak{gl}_2(\mathbb{F})$ via the identification    
\[
 x=\begin{pmatrix}[1.3]
     0 &1~ \\
     0 & 0~\\
 \end{pmatrix},~
 y=\begin{pmatrix}[1.3]
     0 & 0~\\
     1& 0~\\
 \end{pmatrix},~
 z=
 \begin{pmatrix}[1.3]
     1 & 0\\
     0 & -1~
 \end{pmatrix},~
 w=\begin{pmatrix}[1.3]
     1 & 0~\\
     0 & 1~
 \end{pmatrix}.
 \]
 It is easy to see that $[p]_3$ is exactly the natural $p$-map $h\mapsto h^{p}$ on $\mathfrak{gl}_{2}(\mathbb{F})$.
\end{remark}

\begin{lemma}\label{adT-2}
For any  $a, b, c\in \mathbb{F}$, one has in $\mathfrak{gl}_2[T]$ that
\begin{enumerate}
    \item $\big(\mad(axT+by)\big)^{p-1}(ax) = \frac{1}{2}(ab)^{\frac{p-1}{2}}\left(axT^{\frac{p-1}{2}}-byT^{\frac{p-3}{2}}\right)$;
    \item $\big (\mad (czT+ax+by)\big)^{p-1}(cz)=ab\left(ab+c^{2}T^{2}\right)^{\frac{p-3}{2}}cz-c^{2}T(ab+c^{2}T^{2})^{\frac{p-3}{2}}(ax+by)$.
\end{enumerate}
\end{lemma}

\begin{proof}
 Note that $\mathbb{F}$ is of characteristic $p$, so $2^{p-1}=1$ in $\mathbb{F}$. Therefore, it suffices to show   
 \begin{enumerate}
     \item $\big(\mad(axT+by)\big)^{k}(ax)=2^{k-1}(ab)^{\frac{k}{2}}\left(axT^{\frac{k}{2}}-byT^{\frac{k-2}{2}}\right)$ and
     \item $\big(\mad (czT+ax+by)\big)^{k}(cz)=2^{k}ab(ab+c^{2}T^{2})^{\frac{k-2}{2}}cz-2^{k}c^{2}T(ab+c^{2}T^{2})^{\frac{k-2}{2}}(ax+by)$
 \end{enumerate}
 for even integer $k> 0$. We prove them by induction on $k$.  When $k=2$, one has
  \begin{equation*}
        [axT+by,[axT+by,ax]]
        =[abz,axT+by]=2ab(axT-by)
  \end{equation*}
  and 
  \begin{align*}
        [czT+ax+by,[czT+ax+by,cz]]&=[czT+ax+by,-2c(ax-by)] \\
        &=4abcz-4c^{2}T(ax+by).
    \end{align*}
Then by the induction hypothesis, one has
    \begin{align*}
        \big(\mad(axT+by)\big)^{k}(ax)
        &=[axT+by,[axT+by,\mad(axT+by))^{k-2}(ax)]]\\
        &=[axT+by,-2^{k-2}(ab)^{\frac{k}{2}}zT^{\frac{k-2}{2}}]\\
        &=2^{k-1}(ab)^{\frac{k}{2}}\left(axT^{\frac{k}{2}}-byT^{\frac{k-2}{2}}\right)
    \end{align*}
    and 
    \begin{align*}
        \big(\mad (czT+ax+by)\big)^{k}(cz) =&[czT+ax+by,[czT+ax+by,\mad (czT+ax+by))^{k-2}(cz)]]\\
        =&[czT+ax+by,-2^{k-1}c(ab+c^{2}T^{2})^{\frac{k-2}{2}}(ax-by)]\\
        =&2^{k}ab(ab+c^{2}T^{2})^{\frac{k-2}{2}}cz-2^{k}c^{2}T(ab+c^{2}T^{2})^{\frac{k-2}{2}}(ax+by)
    \end{align*}
    as desired.
\end{proof}

\begin{lemma}\label{axbyczp}
    Let $[p]$ be a $p$-semilinear map on $\mathfrak{gl}_2$. Then for any $a, b, c, d\in \mathbb{F}$,
    \begin{align}\label{gl2p}  
    (ax+by+cz+dw)^{[p]}=& ~ a^{p}x^{[p]}+b^{p}y^{[p]}+c^{p}z^{[p]} + d^p w^{[p]}\notag\\&+(c^{2}+ab)^{\frac{p-1}{2}}(ax+by)+((c^{2}+ab)^{\frac{p-1}{2}}-c^{p-1})cz.
    \end{align}
\end{lemma}
\begin{proof} Since $w$ belongs to the center of $\mathfrak{gl}_2$, it suffices to show   \begin{align} 
    (ax+by+cz)^{[p]}=& ~ a^{p}x^{[p]}+b^{p}y^{[p]}+c^{p}z^{[p]} \notag\\&+(c^{2}+ab)^{\frac{p-1}{2}}(ax+by)+\left((c^{2}+ab)^{\frac{p-1}{2}}-c^{p-1}\right)cz. \notag
    \end{align}   
    By Lemma \ref{adT-2} (2) and the equality \eqref{adT-1}, one has
    \begin{align*}
        (ax+by+cz)^{[p]}
        &=(ax+by)^{[p]}+c^{p}z^{[p]}
        +ab\left(\sum_{n=0}^{\frac{p-3}{2}}\frac{\tbinom{\frac{p-3}{2}}{n}}{2n+1}(c^{2})^{n}(ab)^{\frac{p-3}{2}-n}\right)cz\\
        & \quad  +c^{2}\left(-\sum_{n=0}^{\frac{p-3}{2}}\frac{\tbinom{\frac{p-3}{2}}{n}}{2n+2}(c^{2})^{n}(ab)^{\frac{p-3}{2}-n}\right)(ax+by).
    \end{align*}
    Note that in $\mathbb{F}$ one has
    \[
    \frac{\tbinom{\frac{p-3}{2}}{n}}{2n+1}=\frac{\frac{p-1}{2}-n}{\frac{p-1}{2}(2n+1)}\tbinom{\frac{p-1}{2}}{n}=\tbinom{\frac{p-1}{2}}{n} \quad \text{and} \quad \frac{\tbinom{\frac{p-3}{2}}{n}}{2n+2}=\frac{1}{2(\frac{p-1}{2})}\tbinom{\frac{p-1}{2}}{n+1}=-\tbinom{\frac{p-1}{2}}{n+1}.
    \]
    So
    \begin{align*}
    ab\left(\sum_{n=0}^{\frac{p-3}{2}}\frac{\tbinom{\frac{p-3}{2}}{n}}{2n+1}(c^{2})^{n}(ab)^{\frac{p-3}{2}-n}\right)
    &=\sum_{n=0}^{\frac{p-3}{2}}\tbinom{\frac{p-1}{2}}{n}(c^{2})^{n}(ab)^{\frac{p-1}{2}-n}\\
    &=(c^{2}+ab)^{\frac{p-1}{2}}-c^{p-1}
    \end{align*}
    and
    \begin{align*}
        c^{2}\left(-\sum_{n=0}^{\frac{p-3}{2}}\frac{\tbinom{\frac{p-3}{2}}{n}}{2n+2}(c^{2})^{n}(ab)^{\frac{p-3}{2}-n}\right)
        &=-\sum_{n=0}^{\frac{p-3}{2}}\left(-\tbinom{\frac{p-1}{2}}{n+1}\right)(c^{2})^{n+1}(ab)^{\frac{p-3}{2}-n} \\
        &=(c^{2}+ab)^{\frac{p-1}{2}}-(ab)^{\frac{p-1}{2}}.
    \end{align*}
   Therefore,
    \begin{align*}
        (ax+by+cz)^{[p]}
        =&(ax+by)^{[p]}+c^{p}z^{[p]}+\left((c^{2}+ab)^{\frac{p-1}{2}}-(ab)^{\frac{p-1}{2}}\right)(ax+by)\\
        &+\left((c^{2}+ab)^{\frac{p-1}{2}}-c^{p-1}\right)cz. 
    \end{align*}
    In addition, by Lemma \ref{adT-2} (1) and the equality \eqref{adT-1}, one has
    \begin{align}\label{gl2p-2}
        (ax+by)^{[p]}&=a^{p}x^{[p]}+b^{p}y^{[p]}+\frac{1}{2}(ab)^{\frac{p-1}{2}}\left(\frac{1}{1+\frac{p-1}{2}}ax-\frac{1}{1+\frac{p-3}{2}}by\right)  \notag \\
        &=a^{p}x^{[p]}+b^{p}y^{[p]}+(ab)^{\frac{p-1}{2}}(ax+by).
    \end{align}
   The desired formula follows immediately.
\end{proof}

\begin{lemma}\label{varphixyzp}
    Let $\varphi:\mathfrak{gl}_2\to \mathfrak{gl}_2$ be an isomorphism of Lie algebras given by the assignments as in \eqref{AutoLk}. Then for any $p$-semilinear map $[p]$ on $\mathfrak{gl}_2$, one has
    \begin{align}
    \varphi(x)^{[p]}&=a_{1}^{p}x^{[p]}+a_{2}^{p}y^{[p]}+a_{3}^{p}z^{[p]}-a_{3}^{p}z,\label{phixp}\\
    \varphi(y)^{[p]}&=b_{1}^{p}x^{[p]}+b_{2}^{p}y^{[p]}+b_{3}^{p}z^{[p]}-b_{3}^{p}z,\label{phiyp}\\
    \varphi(z)^{[p]}&=c_{1}^{p}x^{[p]}+c_{2}^{p}y^{[p]}+c_{3}^{p}z^{[p]}+c_{1}x+c_{2}y+(c_{3}-c_{3}^{p})z.\label{phizp}
    \end{align}
\end{lemma}

\begin{proof}
    Let 
    \[
    A=\begin{pmatrix}[1.3]
        a_{1} & b_{1} & c_{1} \\
        a_{2}  & b_{2} & c_{2} \\
        a_{3} & b_{3}  & c_{3}
    \end{pmatrix}.
    \]
    Then by equations \eqref{eqisogl2}, the adjoint matrix of A is 
    \[
    A^{*}=
    \begin{pmatrix}[1.3]
        b_{2} & b_{1} &2b_{3} \\
        a_{2}  & a_{1} & 2a_{3}\\
        \frac{1}{2}c_{2}  &\frac{1}{2}c_{1}  &c_{3}
    \end{pmatrix}.
    \]
    By a direct calculation, one has $\det A^{*}=\det A$. So
    \[
    (\det A)^{2}=\det(A)\det(A^*)=\det(AA^{*})=(\det A)^{3}.
    \]
    Since $\varphi$ is invertible, one obtains $\det A=1$. Therefore $A A^{*}=I_{3}$, which follows that
    \begin{eqnarray}\label{more-equations}
    a_{1}a_{2} +a_{3}^{2}=0, \quad  b_{1}b_{2}+b_{3}^{2}=0,\quad c_{1}c_{2}+c_{3}^{2}-1 =0. 
    \end{eqnarray}
    Now by Lemma \ref{axbyczp},
    \begin{align*}
\varphi(z)^{[p]}&=c_{1}^{p}x^{[p]}+c_{2}^{p}y^{[p]}+c_{3}^{p}z^{[p]}+(c_{3}^{2}+c_{1}c_{2})^{\frac{p-1}{2}}(c_{1}x+c_{2}y) \\
  & \quad +\left((c_{3}^{2}+c_{1}c_{2})^{\frac{p-1}{2}}-c_{3}^{p-1}\right)c_{3}z\\
&=c_{1}^{p}x^{[p]}+c_{2}^{p}y^{[p]}+c_{3}^{p}z^{[p]}+c_{1}x+c_{2}y+(c_{3}-c_{3}^{p})z.
    \end{align*}
    The other two equations can be derived via similar computations.
\end{proof}

\begin{remark}\label{Groebner}
  We would like to illustrate an algorithmic   way to obtain the equations in \eqref{more-equations}. It may be useful in similar researches. Let $I_{\mathbb{Q}}$ be the ideal of $R_{\mathbb{Q}}=\mathbb{Q}[t,a_{1},a_{2},a_{3},b_{1},b_{2},b_{3},c_{1},c_{2},c_{3}]$
    generated by the left-hand sides of those equations in \eqref{eqisogl2} together with
    \[
    t(a_{1}b_{2}c_{3}+a_{2}b_{3}c_{1}+a_{3}b_{1}c_{2}-a_{1}b_{3}c_{2}-a_{2}b_{1}c_{3}-a_{3}b_{2}c_{1})-1.
    \]
    Using Magma, we obtain a Gr\"obner basis of $I_{\mathbb{Q}}$ with respect to the lexicographical order  determined by $t>a_{1}>a_{2}>a_{3}>b_{1}>b_{2}>b_{3}>c_{1}>c_{2}>c_{3}.$ See Listing  \ref{autoN} for the codes.  It follows that
    \begin{equation*}
       a_{1}a_{2} +a_{3}^{2}, \quad  b_{1}b_{2}+b_{3}^{2},\quad c_{1}c_{2}+c_{3}^{2}-1 \quad \in \quad I_{\mathbb{Q}},
    \end{equation*}
    because these three relations are all contained in the resulting Gr\"obner basis.
    Though the ground field is $\mathbb{Q}$, the conclusion  holds if we replace it by $\mathbb{F}$. To see this, let us denote by $f_i$ the input polynomials. For any element $g$ of the resulting Gr\"obner basis, there exist polynomials $h_{i}\in R_{\mathbb{Q}}$ such that $g =\sum h_{i}f_{i}.$ One may find $h_i$ by using the Magma Function Coordinates. If the coefficients of $h_i$ are all definable over $\mathbb{F}$, then one may view $g =\sum h_{i}f_{i}$ as an equality in $R_{\mathbb{F}}:=\mathbb{F}[t,a_{1},a_{2},a_{3},b_{1},b_{2},b_{3},c_{1},c_{2},c_{3}]$. For example, it turns out in $R_{\mathbb{Q}}$ that 
    \begin{equation*}
        a_{3}^{2}+a_{1}a_{2}=-\frac{1}{2}a_{3}(a_{2}c_{1}-a_{1}c_{2}-2a_{3})-\frac{1}{2} a_{2}(a_{1}c_{3}-a_{3}c_{1}-a_{1})- \frac{1}{2}a_{1}(a_{3}c_{2} - a_{2}c_{3}-a_{2}).
    \end{equation*}
    Since $\mathbb{F}$ is of characteristic $\geq 3$, it follows that $\frac{1}{2}$ is definable over $\mathbb{F}$ and hence this equality also holds in $R_{\mathbb{F}}$. The verifications for $b_{1}b_{2}+b_{3}^{2}$ and $c_{1}c_{2}+c_{3}^{2}-1$ are similar.
\end{remark}

\begin{lemma}\label{iso(3)and(1)or(8)}
    Let $\mu, \lambda\in \mathbb{F}$ and let $[p]$ be the $p$-map on $\mathfrak{gl}_2$ given by
    \begin{equation*}
    x^{[p]}=w,\quad y^{[p]}=\mu w, \quad z^{[p]}=z+w,\quad w^{[p]}=\lambda w.
    \end{equation*}
    \begin{enumerate}
       \item If $\mu=-\frac{1}{4}$ and $\lambda=0$, then $[p]$ is conjugate to $[p]_2$;
        \item If $\mu=-\frac{1}{4}$ and $\lambda\neq0$, then $[p]$ is conjugate to  $[p]_4$;
        \item If $\mu\neq -\frac{1}{4}$, then $[p]$ is conjugate to $[p]_{5, \lambda'}$ with $\lambda'=(4\mu+1)^{\frac{p-1}{2}}\lambda$.
    \end{enumerate}
\end{lemma}
\begin{proof}
Firstly, assume $\mu\neq -\frac{1}{4}$. Let $\varphi: \mathfrak{gl}_2 \to \mathfrak{gl}_2$ be the linear map given by
    \begin{align*}
        x&\mapsto x-\frac{1}{s^{2}}y+\frac{1}{s}z &
        y&\mapsto -\frac{(s-1)^{2}}{4}x+\frac{(s+1)^{2}}{4s^{2}}y+\frac{s^{2}-1}{4s}z\\
        z&\mapsto (1-s)x-\frac{s+1}{s^{2}}y+\frac{1}{s}z &
        w&\mapsto \frac{1}{s^{p}}w,
    \end{align*}
   where $s=(4\mu+1)^{\frac{1}{2p}}$. By checking the equations in \eqref{eqisogl2}, it is straightforward to see $\varphi \in \mAut(\mathfrak{gl}_2)$.  By Lemma \ref{varphixyzp}, one may readily verify that $\varphi(x)^{[p]_{5,\lambda'}}=\varphi(x^{[p]})$, $\varphi(y)^{[p]_{5,\lambda'}}=\varphi(y^{[p]})$, $\varphi(z)^{[p]_{5,\lambda'}}=\varphi(z^{[p]})$ and $\varphi(w)^{[p]_{5,\lambda'}}=\varphi(w^{[p]})$. For example, one has
    \begin{align*}
        \varphi(z)^{[p]_{5,\lambda'}}
        &=\frac{1}{s^{p}}(z+w)+\left(\frac{1}{s}-\frac{1}{s^{p}}\right)z+(1-s)x-\frac{s+1}{s^{2}}y\\
        &=(1-s)x-\frac{s+1}{s^{2}}y+\frac{1}{s}z+\frac{1}{s^{p}}w\\
        &=\varphi(z+w)=\varphi(z^{[p]}).
    \end{align*}
    Thus $\varphi\circ[p]\circ\varphi^{-1}=[p]_{5,\lambda'}$, that is, $[p]$ is conjugate to $[p]_{5,\lambda'}$. 
  
  Next assume $\mu=-\frac{1}{4}$ and $\lambda=0$ (resp. $\lambda\neq0$). Let $\psi:\mathfrak{gl}_2\to \mathfrak{gl}_2$ be the linear map given by
      \begin{align*}
          & x\mapsto -4x+y-2z,\quad y\mapsto x,\quad z\mapsto -z-4x,\quad w\mapsto -4w.\\
        \big(\text{resp.} \quad & x\mapsto -4sx+ \frac{1}{s}y -2z, \quad y\mapsto sx, \quad 
        z\mapsto -4sx - z, \quad    w\mapsto -4s^p w, \big)
      \end{align*}
      where $s=\lambda^{\frac{1}{p(p-1)}}$. It is also straightforward to show as above that $\psi \in \mAut(\mathfrak{gl}_2)$  and $ \psi\circ [p]\circ\psi^{-1} = [p]_{2}$ (resp. $ \psi\circ [p]\circ\psi^{-1} = [p]_{4}$). Therefore, $[p]$ is conjugate to $[p]_{2}$ (resp. $[p]_4$).
\end{proof}

\begin{remark}\label{Groebner-2}
   Indeed, Magma is used to find the critical number $\mu=- \frac{1}{4}$ in Lemma \ref{iso(3)and(1)or(8)} as well as the isomorphisms $\varphi$ and $\psi$ in the proof. At the very beginning, we guess that $\varphi\circ[p]=[p]_{\lambda'}\circ\varphi$ for some $\varphi\in \mAut(\mathfrak{gl}_{2})$ and some  $\lambda'\in \mathbb{F}$.  
   Describe $\varphi(x), \varphi(y), \varphi(z), \varphi(w)$ as that in \eqref{AutoLk} and then calculate the value of $\varphi\circ[p]$ and $[p]_{\lambda'}\circ\varphi$ on $x$, $y$, $z$,  $w$, one gets 
    \begin{equation}\label{constructofphi}
    d_{4}-a_{3}^{p}=0,\quad \mu d_{4}-b_{3}^{p}=0,\quad d_{4}-c_{3}^{p}=0,\quad d_{4}^{p-1}\lambda'-\lambda=0.
    \end{equation}
    Let $J_{\mathbb{Q}}$ be the ideal of $S_{\mathbb{Q}}=\mathbb{Q}[t,a_{1},a_{2},a_{3},b_{1},b_{2},b_{3},c_{1},c_{2},c_{3},X, Y]$ generated by the left-hand sides of those equations in \eqref{eqisogl2} together with $a_{3}-X$, $b_{3}-XY$, $c_{3}-X$ and
    \[
    t(a_{1}b_{2}c_{3}+a_{2}b_{3}c_{1}+a_{3}b_{1}c_{2}-a_{1}b_{3}c_{2}-a_{2}b_{1}c_{3}-a_{3}b_{2}c_{1})X-1.
    \]
    Here, we employ $X$ and $Y$ to deal with $d_{4}^{\frac{1}{p}}$ and $\mu^{\frac{1}{p}}$ respectively, which cannot be recognized as constants by Magma program because $p$ is not determined. Using Magma, we obtain a Gr\"obner basis of $J_{\mathbb{Q}}$ with respect to the lexicographical order determined by $t>a_{1}>a_{2}>a_{3}>b_{1}>b_{2}>b_{3}>c_{1}>c_{2}>c_{3}>X>Y$.
    See Listing \ref{autoresofN} for the codes. As a result of the output, one has
    \begin{eqnarray*}
    X^2Y +\frac{1}{4}X^2 - \frac{1}{4} \quad \in \quad  J_{\mathbb{Q}}.
    \end{eqnarray*}
    It is not hard to check that the conclusion also holds if the ground field $\mathbb{Q}$ is replaced by $\mathbb{F}$, as explained in Remark \ref{Groebner}.
    Let $X=d_{4}^{\frac{1}{p}}$ and $Y=\mu^{\frac{1}{p}}$, one obtains that $d_4$ and $\mu$ need to satisfy $d_{4}^{\frac{2}{p}}\mu^{\frac{1}{p}} + \frac{1}{4}d_{4}^{\frac{2}{p}} - \frac{1}{4}=0$.
    Since $4^p=4$ in $\mathbb{F}$, it follows that $$d_{4}^{2}(4\mu+1)=1.$$ Therefore, we need to assume $\mu\neq-\frac{1}{4}$ for our purpose. Then take $d_{4}=(4\mu+1)^{-\frac{1}{2}}$,  $a_{1}=1$ and put them in the equations in \eqref{eqisogl2} and \eqref{constructofphi}, one gets the value of $a_i, b_i, c_i$ and $\lambda'$ immediately. The above argument illustrates how $\varphi$ has been constructed. Similar discussion for $\psi$.
\end{remark}

Now we are ready to prove Proposition \ref{Lkp}.
\begin{proof}[Proof of Proposition \ref{Lkp}]
Let $[p]$ be a $p$-map on $\mathfrak{gl}_2$. As explained in the paragraph before Proposition \ref{Lkp}, we may assume
$x^{[p]}=\alpha w,  y^{[p]}=\beta w, z^{[p]}=z+\gamma w, w^{[p]}=\delta w $ for some $\alpha, \beta, \gamma, \delta\in \mathbb{F}$. The discussion is broke down into following situations up on $\alpha$, $\beta$, $\gamma$ and $\delta$. 
\begin{enumerate}
    \item Assume $\alpha=\beta=\gamma=0$ and $\delta=0$. Then $[p]=[p]_{1}$.
    
    \item Assume $\alpha=\beta=\gamma=0$ and $\delta\neq0$. Then $\sigma_{1,1,1,d}\circ [p]\circ \sigma_{1,1,1,d}^{-1}=[p]_{3}$, where $d=\delta^{\frac{1}{p-1}}$.
    
    \item Assume $\alpha=\beta=0$ and $\gamma\neq 0$. Then $\sigma_{1,1,1,\gamma^{-1}}\circ [p]\circ \sigma_{1,1,1,\gamma^{-1}}^{-1}=[p]_{5,\gamma^{p-1}\delta}$.

    \item Assume $\alpha\neq0$, $\beta=\gamma=0$ and $\delta=0$. Then $\sigma_{1,1,1,\alpha^{-1}}\circ [p]\circ \sigma_{1,1,1,\alpha^{-1}}^{-1}=[p]_{2}$.
    
    \item Assume $\alpha\neq0$, $\beta=\gamma=0$ and $\delta \neq 0$. Then $\sigma_{a,b,1,d}\circ [p]\circ \sigma_{a,b,1,d}^{-1}=[p]_{4}$, where  $a=\left(\delta^{\frac{1}{p-1}}\alpha\right)^{\frac{1}{p}}$, $b=\left(\delta^{\frac{1}{p-1}} \alpha\right)^{-\frac{1}{p}}$, $d=\delta^{\frac{1}{p-1}}$.
    
   \item Assume $\alpha \neq 0$, $\beta=0$ and $ \gamma\neq 0$. Then $\varphi\circ[p]\circ\varphi^{-1}=[p]_{5,\gamma^{p-1}\delta}$, where $\varphi:\mathfrak{gl}_2\to \mathfrak{gl}_2$ is the isomorphism of Lie algebras  given by
\[
x\mapsto \left(\frac{\alpha}{\gamma}\right)^{\frac{1}{p}}(-x+y+z), \quad y\mapsto -\left(\frac{\gamma}{\alpha}\right)^{\frac{1}{p}}y, \quad z\mapsto 2y+z, \quad w\mapsto \gamma^{-1}w.
\]

\item Assume $\alpha=0$ and $\beta\neq0$. Then $[p]':=\Delta\circ[p]\circ\Delta^{-1}$ satisfies that
\[
x^{[p]'}=\beta w,\quad y^{[p]'}=\alpha w=0,\quad z^{[p]'}=z+\gamma w,\quad w^{[p]'}=\delta w.
\]  
Here, $\Delta:\mathfrak{gl}_2\to \mathfrak{gl}_2$ is the isomorphism of Lie algebras given by
\[
x\mapsto -y,\quad y\mapsto -x,\quad z\mapsto -z,\quad w\mapsto -w.
\]
Clearly, $[p]'$  falls into the situations (4), (5) or (6).

\item Assume $\alpha\neq 0$, $\beta\neq 0$ and $\gamma=0$. Then  $\psi\circ[p]\circ\psi^{-1}=[p]_{5,s^{p-1}\delta}$, where $s= \sqrt{\alpha\beta}$ and $\psi:\mathfrak{gl}_2\to \mathfrak{gl}_2$ is the isomorphism of Lie algebras given by
\[
x\mapsto -\frac{1}{2}\left(\frac{\alpha}{s}\right)^{\frac{1}{p}}(x-y-z), \quad y\mapsto \frac{1}{2}\left(\frac{\beta}{s}\right)^{\frac{1}{p}}(x-y+z), \quad z\mapsto x+y, \quad w\mapsto \frac{2}{s}w.
\]

\item Assume $\alpha,\beta,\gamma\neq0$. Let  $[p]':=\sigma_{a, b, 1,d}\circ[p]\circ\sigma_{a, b, 1 ,d}^{-1}$, where $a=\left(\frac{\alpha}{\gamma}\right)^{\frac{1}{p}}$, $b=\left(\frac{\gamma}{\alpha}\right)^{\frac{1}{p}}$ and $d=\gamma^{-1}$. Then $[p]'$  satisfies that
\begin{equation*}
    x^{[p]'}=w,\quad y^{[p]'}=\beta w, \quad z^{[p]'}=z+w,\quad w^{[p]'}=\delta w.
\end{equation*} 
By Lemma \ref{iso(3)and(1)or(8)}, $[p]'$ is conjugate to $[p]_{2}$, $[p]_{4}$ or $[p]_{5,\lambda'}$ with $\lambda'=(4\beta+1)^{\frac{p-1}{2}}\delta$.
\end{enumerate}
As a summary, $[p]$ is  conjugate to one of $[p]_1, [p]_2, [p]_3, [p]_4, [p]_{5,\lambda}$, where $\lambda\in \mathbb{F}$.

Next, we proceed to show  the listed $p$-maps are  not conjugate to each other. By a straightforward but tedious calculation on $\dim Z(\mathfrak{gl}_2)^{[p]}$ and $\dim [\mathfrak{gl}_2,\mathfrak{gl}_2]^{[p]}$, we have Table \ref{Lkpdim}. 
\begin{table}[htbp]
\centering
\caption{2-tuple of invariant of $L=\mathfrak{gl}_2$}\label{Lkpdim}
\setlength{\tabcolsep}{5mm}{
\renewcommand\arraystretch{1.4}
\begin{tabular}{c|cc|c}
\hline
$[p]$      & $\dim Z(L)^{[p]}$ & $\dim [L,L]^{[p]}$  & Conditions  \\ \hline
$[p]_{1}$  & 0 & 3 &\\ \hline
$[p]_{2}$ & 0 & 4  &  \\ \hline
$[p]_{3}$ & 1 & 3  &           \\ \hline
$[p]_{4}$ & 1 & 4  &\\\hline
\multirow{2}{*}{$[p]_{5,\lambda}$} & 1 & 4 & $\lambda\neq0$\\\cline{2-4}
       & 0 & 4 &$\lambda=0$ \\\hline
\end{tabular}}
\end{table}
So by Lemma \ref{ConjInvar}, it remains to show that $[p]_{2}$ and $[p]_{5,0}$ are not conjugate; $[p]_{4}$ and $[p]_{5,\lambda}$ are not conjugate for $\lambda\neq 0$; and $[p]_{5,\lambda}$ and $[p]_{5,\mu}$ are not conjugate for $\lambda, \mu\neq 0$ with $\lambda\neq \mu$.

Suppose $[p]_{2}$ is conjugate to $[p]_{5,0}$ (resp. $[p]_{4}$ is conjugate to $[p]_{5,\lambda}$, where $\lambda\neq 0$). Then there exists an isomorphism $\Phi\in \mAut(\mathfrak{gl}_2)$ such that $\Phi\circ[p]_{2}=[p]_{5,0}\circ\Phi$ (resp. $\Phi\circ[p]_{4}=[p]_{5,\lambda}\circ\Phi$). Describe $\Phi(x), \Phi(y), \Phi(z), \Phi(w)$ as that in \eqref{AutoLk}.  By Lemma \ref{varphixyzp}, 
\begin{align*}
    \Phi(z^{[p]_{2}})=\Phi(z)^{[p]_{5,0}} \quad (\text{resp. }\Phi(z^{[p]_{4}})=\Phi(z)^{[p]_{5,\lambda}}) & \quad \Longrightarrow \quad 0=c_{3}^{p}w,\\
    \Phi(y^{[p]_{2}})=\Phi (y)^{[p]_{5,0}} \quad (\text{resp. }\Phi(z^{[p]_{4}})=\Phi(z)^{[p]_{5,\lambda}}) & \quad\Longrightarrow \quad 0=b_{3}^{p}w.
\end{align*}
It follows that $c_{3}=b_{3}=0$. Then taking them into the equations in \eqref{eqisogl2}, one gets $b_{1}=b_{2}=0$. But it is impossible because $\Phi$ is invertible. 

Suppose $[p]_{5,\lambda}$ is conjugate to $[p]_{5,\mu}$, where $\lambda, \mu\neq 0$. Then there is $\Psi\in \mAut(\mathfrak{gl}_2)$  such that $\Psi\circ[p]_{5,\lambda}=[p]_{5,\mu}\circ\Psi$. Describe $\Psi(x), \Psi(y), \Psi(z), \Psi(w)$  as that in \eqref{AutoLk}.  By Lemma \ref{varphixyzp},
\begin{align*}
    \Psi(x^{[p]_{5,\lambda}})=\Psi(x)^{[p]_{5,\mu}} &\quad \Longrightarrow \quad 0=a_{3}^{p}w,\\
     \Psi(y^{[p]_{5,\lambda}})=\Psi(y)^{[p]_{5,\mu}}&\quad \Longrightarrow \quad 0=b_{3}^{p}w,\\
     \Psi(z^{[p]_{5,\lambda}})=\Psi(z)^{[p]_{5,\mu}} &\quad \Longrightarrow \quad d_{4}w=c_{3}^{p}w,\\
    \Psi(w^{[p]_{5,\lambda}})=\Psi(w)^{[p]_{5,\mu}}    
    &\quad \Longrightarrow \quad d_{4}\lambda w=d_{4}^{p}\mu w.
\end{align*}
So $a_{3}=b_{3}=0$, $d_{4}=c_{3}^{p}$ and $\lambda=d_{4}^{p-1}\mu $. Taking $a_{3}=b_{3}=0$ into the equations in \eqref{eqisogl2}, it follows that $c_{3}^{2}=1$ and hence $d_{4}^{p-1}=1$. Therefore, $\lambda=\mu$.
\end{proof}

\section{The classification theorem}\label{classification-thm}

We are now in position to state and prove the main result of this paper. Recall that the base field $\mathbb{F}$ is assumed to be algebraically closed of characteristic $p>0$.

\begin{thm}\label{theclassification}
Every restricted Lie algebra of dimension $4$  is isomorphic to one of the examples that listed  in Table \ref{mainresult} (see Appendix \ref{tableres}). Moreover, these restricted Lie algebras that listed  in  Table \ref{mainresult}  are not isomorphic to each other except that
\begin{enumerate}
    \item $L_{2}^{15}(\lambda_{1})\cong L_{2}^{15}(\lambda_{2})$ if and only if $\lambda_{1}^{p(p-1)}(\lambda_{2}^{p-1}+1)^{p+1} =\lambda_{2}^{p(p-1)}(\lambda_{1}^{p-1}+1)^{p+1}$ for any scalars $\lambda_{1}, \lambda_{2}\in \mathbb{F}$;
    \item $L_{4}^{11}(\lambda_{1})\cong L_{4}^{11}(\lambda_{1})$ if and only if  $\left(\begin{smallmatrix} 1 \\ \lambda_{1} \end{smallmatrix}\right)=aA\!\left(\begin{smallmatrix} 1 \\ \lambda_{2} \end{smallmatrix}\right)$  for some nonzero element $a\in \mathbb{F}^\times$ and invertible matrix $A\in \mathrm{GL}_{2}(\mathbb{F}_{p})$, for any scalars $\lambda_1,\lambda_2\in \mathbb{F}$;
    \item $L_{6}^{1}(\xi_{1},\eta_{1}) \cong L_{6}^{1}(\xi_{2},\eta_{2})$ if and only if $(\xi_{1},\eta_{1}) =\tau\cdot (\xi_{2},\eta_{2})$ for some permutation $\tau \in S_{3}$, for any scalars $\xi_{1}, \eta_{1}, \xi_{2}, \eta_{2} \in \mathbb{F}_{p}^{\times}$.
\end{enumerate}
\end{thm}

\begin{proof}
In Section \ref{Classify-Lie}, we obtain a complete representatives of Lie algebras of dimension $4$. Then in Sections \ref{Classify-restrictedLie-1} - \ref{Classify-restrictedLie-3}, we obtain a complete conjugate representatives of $p$-maps on each of the representative Lie algebras obtained in  Section \ref{Classify-Lie}. Combine all the propositions in Sections \ref{Classify-Lie} - \ref{Classify-restrictedLie-3} as well as Remarks \ref{L5NONiso} and  \ref{L6NONiso} together, one may readily get the  classification theorem.
\end{proof}

\begin{remark}\label{Hopf-group-scheme}
It is well-known that over a field of positive characteristic, the assignments $$(L, [p]) \mapsto \mathfrak{u}(L,[p])$$ define an equivalence from the category of restricted Lie algebras to that of primitively generated connected Hopf algebras, and the assignments $$(L, [p]) \mapsto {\rm Spec} (\mathfrak{u}(L,[p])^{*})$$ define an equivalence from the category of finite dimensional  restricted Lie algebras to that of infinitesimal group schemes of height at most $1$. Here, $\mathfrak{u}(L,[p])$ denotes the restricted universal enveloping algebra of $(L, [p])$, and $\mathfrak{u}(L,[p])^{*}$ the dual Hopf algebra of $\mathfrak{u}(L,[p])$.  Since $$\dim \mathfrak{u}(L,[p]) = p^{\dim L},$$
Theorem \ref{theclassification} provides a classification of primitively generated connected Hopf algebras over $\mathbb{F}$ of dimension $p^4$ and  of  infinitesimal group schemes  over $\mathbb{F}$ of order $p^{4}$ and height $1$.  
\end{remark}

Note that the isomorphism classes of the classification given in Theorem \ref{classification-thm} 
consists of $2$ (resp. $3$) infinite one-parameter families, and a finite number of individuals or finite parameterized families for $p=2$ (resp. $p\geq 3$).  These infinite families are represented by $L_{2}^{15}(\lambda)$ and $L_{4}^{11}(\lambda)$ (resp. $L_{2}^{15}(\lambda)$, $L_{4}^{11}(\lambda)$ and $\mathfrak{gl}_{2}^{5}(\lambda)$) with $\lambda$ runs over $\mathbb{F}$, as listed  in  Table \ref{mainresult}.

In the remaining of this section, we proceed  to describe the number of isomorphism classes that do not belong to the infinite one-parameter families. To this end, recall from \eqref{actionS3}  that the action of $S_{3}$ on $\mathbb{F}_{p}^{\times} \times \mathbb{F}_{p}^{\times}$ is given by $ (12)\cdot(\xi,\eta)=(\xi^{-1},\xi^{-1}\eta)$ and $(23)\cdot(\xi,\eta)=(\eta,\xi)$. The next result counts the number of orbits of this action.

\begin{lemma}\label{S3orbit}
    Let $\Omega$ be the set of all orbits of the action of $S_{3}$ on $\mathbb{F}_{p}^{\times} \times \mathbb{F}_{p}^{\times}$. Then
    \begin{equation*}
        \#\Omega=
        \begin{cases}
            \frac{p^{2}+p}{6}, &\quad \text{when $3\nmid (p-1)$}\\
            \frac{p^{2}+p+4}{6}, &\quad \text{when $3\mid (p-1)$}.
        \end{cases}
    \end{equation*}
\end{lemma}
\begin{proof}
    Let $X=\mathbb{F}_{p}^{\times} \times \mathbb{F}_{p}^{\times}$ and $X^{\sigma}=\{(\xi,\eta)\in X~|~\sigma\cdot (\xi,\eta)=(\xi,\eta)\}$ for $\sigma\in S_{3}$. It is easy to see that $\#X^{\mathrm{id}}=\#X=(p-1)^{2}$. Note that 
    \[
    (12)(\xi,\eta)=(\xi^{-1},\xi^{-1}\eta),\quad (23)(\xi,\eta)=(\eta,\xi),\quad (13)(\xi,\eta)=(\eta^{-1}\xi,\eta^{-1}).
    \]
    So $\#X^{(12)}=\#X^{(23)}=\#X^{(13)}=\#\mathbb{F}_{p}^{\times}=p-1$. In addition, one has 
    \[
    (123)(\xi,\eta)=(\eta^{-1},\eta^{-1}\xi),\quad (132)(\xi,\eta)=(\xi^{-1}\eta,\xi^{-1}).
    \]
    It follows that $\#X^{(123)} =\#X^{(132)} = \# \{\eta\in \mathbb{F}_{p}^{\times } ~| ~ \eta^{3} =1\}$. By the standard results of finite group theory, $\mathbb{F}_{p}^{\times}$ admits an element of order $3$ if and only if $3\mid (p-1)$. In addition,  if $\alpha \in \mathbb{F}_{p}^{\times}$ is an element of order $3$ then so is $\alpha^{2}$. 
    Therefore one has
    \[
    \# X^{(123)} =\#X^{(132)}=
    \begin{cases}
        1, & \text{when $3\nmid p-1$;}\\
        3, & \text{when $3\mid p-1$.}
    \end{cases}
    \]
    Since $\#\Omega=\frac{\sum_{\sigma\in S_{3}}\#X^{\sigma}}{\# S_{3}}$ by the Burnside's Lemma, the desired formula  follows readily.
\end{proof}

Let $N_{p}$ be the number of isomorphism classes of restricted Lie algebras of dimension $4$ which do not belong to the infinite families that represented by $L_{2}^{15}(\lambda)$ and $L_{4}^{11}(\lambda)$ for $p=2$ (resp. $L_{2}^{15}(\lambda)$, $L_{4}^{11}(\lambda)$ and $\mathfrak{gl}_{2}^{5}(\lambda)$ for $p\geq 3$), where $\lambda$ runs over $\mathbb{F}$, as listed in Tabel \ref{mainresult}.

\begin{coro}\label{numberisoclass}
$N_{2} =42$, $N_{3} =63$, and for $p\geq 5$ one has
\begin{equation*}
        N_{p}=
        \begin{cases}
            \frac{p^{2}+28p+291}{6}, &\quad \text{when $3\nmid (p-1)$};\\
            \frac{p^{2}+28p+295}{6}, &\quad \text{when $3\mid (p-1)$}.
        \end{cases}
    \end{equation*}
\end{coro}

\begin{proof}
By Theorem \ref{theclassification}, it is clear that $N_{2} =42$ and $N_{3} =63$. Now assume $p\geq 5$.  We first  count the number of the (distinct)  isomorphism classes in the finite parameterized families that listed in Table \ref{mainresult}, which are represented by $L_{5}^{i}(\xi)$ for $i=1, \ldots,8$, $L_{6}^{1}(\xi,\eta)$ and $N_{3}^{1}(\xi)$, with $\xi, \eta$ run over certain subsets of $\mathbb{F}_{p}^{\times}$. Note that $\#\Xi_{p} = \frac{p+1}{2}$ and $\#Q_{p} = \frac{p-1}{2}$. 
By Theorem \ref{theclassification} and Lemma \ref{S3orbit}, one may see that $L_{5}^{1}(\xi)$, $L_{5}^{2}(\xi)$, $L_{5}^{5}(\xi)$, and $L_{5}^{7}(\xi)$ exhibit $\frac{p+1}{2}$ isomorphism classes, $L_{5}^{4}(\xi)$ and $L_{5}^{8}(\xi)$ exhibit $\frac{p-1}{2}$ isomorphism classes, $L_{5}^{3}(\xi)$ and $L_{5}^{7}(\xi)$ exhibit $\frac{p-3}{2}$ isomorphism classes, $L_{6}^{1}(\xi, \eta)$ exhibits $\frac{p^{2}+p}{6}$ (resp. $\frac{p^{2}+p+4}{6}$) isomorphism classes when $3\nmid (p-1)$ (resp. $3\mid (p-1)$), and $N_{3}^{1}(\xi)$ exhibits $\frac{p-1}{2}$ isomorphism classes. In addition,  there are $53$ individual isomorphism classes that listed in Table \ref{mainresult}. Add up the above amounts, the result follows.
\end{proof}

\section*{Acknowledgments}

This research is supported by the NSFC (Grant Nos. 12371039).

\appendix

\newpage

\section{Restricted Lie algebras of dimension $4$}\label{tableres}
\setlength\LTleft{-0.9cm}
\begin{ThreePartTable}
  \begin{TableNotes}\footnotesize
\item[$\dagger$] $\Xi_{p}=\{\xi_{p}^{r}~|~0\leq r\leq \frac{p-1}{2}\}$ for a fix generator $\xi_{p}\in \mathbb{F}_{p}^{\times}$;
\item[$\ddagger$] $Q_{p}=\{\alpha\in\mathbb{F}_{p}^{\times}~|~\alpha^{\frac{p-1}{2}}=1\}$.
  \end{TableNotes}
\small
{
\begin{longtable}{c|c|c|c|c|c}
\caption{Restricted Lie algebras of dimension $4$}\label{mainresult} \\

\hline \textbf{Row} & \textbf{Name} & \textbf{Lie structure} & \textbf{$[p]$-structure} & \textbf{Parameters} & \textbf{Char.}  \\ \hline
\endfirsthead

\multicolumn{6}{c}%
{{\tablename\ \thetable{} -- continued from previous page}} \\
\hline \textbf{Row} & \textbf{Name} & \textbf{Lie structure} & \textbf{$[p]$-structure} & \textbf{Parameters} & \textbf{Char.}  \\ \hline
\endhead

\hline \multicolumn{6}{r}{{Continued on next page}} \\ 
\endfoot

\hline \insertTableNotes
\endlastfoot
\hline
1 & $L_{1}^{1}$ & \multirow{12}{*}{Abelian} & Trivial &  &  \\ \cline{1-2} \cline{4-6} 
2 & $L_{1}^{2}$ &  & $x^{[p]}=y$ &  &  \\ \cline{1-2} \cline{4-6} 
3 & $L_{1}^{3}$ &  & $x^{[p]}=y,~ z^{[p]}=w$ &  &  \\ \cline{1-2} \cline{4-6} 
4 & $L_{1}^{4}$ &  & $x^{[p]}=y,~  y^{[p]}=z$ &  &  \\ \cline{1-2} \cline{4-6} 
5 & $L_{1}^{5}$ &  &$x^{[p]}=y,~ y^{[p]}=z,~ z^{[p]}=w$  &  &  \\ \cline{1-2} \cline{4-6} 
6 & $L_{1}^{6}$ &  &$x^{[p]}=x$  &  &  \\ \cline{1-2} \cline{4-6} 
7 & $L_{1}^{7}$ &  &$x^{[p]}=x,~y^{[p]}=z$  &  &  \\ \cline{1-2} \cline{4-6} 
8 & $L_{1}^{8}$ &  &$x^{[p]}=x,~ y^{[p]}=y$  &  &  \\ \cline{1-2} \cline{4-6} 
9 & $L_{1}^{9}$ &  & $x^{[p]}=x,~y^{[p]}=z,~z^{[p]}=w$ &  &  \\ \cline{1-2} \cline{4-6} 
10 & $L_{1}^{10}$ &  &$x^{[p]}=x,~  y^{[p]}=y,~ z^{[p]}=w$  &  &  \\ \cline{1-2} \cline{4-6} 
11 & $L_{1}^{11}$ &  &$x^{[p]}=x,~ y^{[p]}=y,~ z^{[p]}=z$  &  &  \\ \cline{1-2} \cline{4-6} 
12 & $L_{1}^{12}$ &  & \makecell{$x^{[p]}=x,~ y^{[p]}=y,~ z^{[p]}=z$,\\ $~ w^{[p]}=w$} &  &  \\ \hline
13 & $L_{2}^{1}$ & \multirow{11}{*}{$[w,x]=y$} & Trivial &  & $p\geq3$ \\ \cline{1-2} \cline{4-6} 
14 & $L_{2}^{2}$ &  & $w^{[p]}=y$ &  &  \\ \cline{1-2} \cline{4-6} 
15 & $L_{2}^{3}$ &  & $w^{[p]}=z$ &  &  \\ \cline{1-2} \cline{4-6} 
16 & $L_{2}^{4}$ &  &  $x^{[p]}=z,~   w^{[p]}=y$ &  &  \\ \cline{1-2} \cline{4-6}
17 & $L_{2}^{5}$ &  & $y^{[p]}=z$ &  & $p\geq3$ \\ \cline{1-2} \cline{4-6} 
18 & $L_{2}^{6}$ &   & $x^{[p]}=y,~ y^{[p]}=z$ &  \\ \cline{1-2} \cline{4-6} 
19 & $L_{2}^{7}$ &  & $z^{[p]}=y$ &  &  \\ \cline{1-2} \cline{4-6} 
20 & $L_{2}^{8}$ &  & $x^{[p]}=z,~  z^{[p]}=y$ &  &  \\ \cline{1-2} \cline{4-6} 
21 & $L_{2}^{9}$ &  & $y^{[p]}=y$ &  &  \\ \cline{1-2} \cline{4-6} 
22 & $L_{2}^{10}$ &  & $x^{[p]}=z,~ y^{[p]}=y$ &  &  \\ \cline{1-2} \cline{4-6} 
23 & $L_{2}^{11}$ &  & $y^{[p]}=y+z$ &  &  $p\geq3$ \\ \cline{1-2} \cline{4-6} 
24 & $L_{2}^{12}$ & \multirow{4}{*}{$[w,x]=y$} & $x^{[p]}=z,~ y^{[p]}=y+z$ &  &  \\ \cline{1-2} \cline{4-6} 
25 & $L_{2}^{13}$ &   & $z^{[p]}=z$ &  & $p\geq3$  \\ \cline{1-2} \cline{4-6} 
26 & $L_{2}^{14}$ &  & $x^{[p]}=y,~  z^{[p]}=z$ &  &  \\ \cline{1-2} \cline{4-6} 
27 & $L_{2}^{15}(\lambda)$ &  & $y^{[p]}=y+\lambda z,~ z^{[p]}=z$ & $\lambda\in\mathbb{F}$ &  \\ \hline
28 & $L_{3}^{1}$ & \multirow{5}{*}{$[w,x]=y,~[w,y]=z$} & Trivial &  & $p\geq5$ \\ \cline{1-2} \cline{4-6} 
29 & $L_{3}^{2}$ &  & $w^{[p]}=z$ &  & $p\geq3$ \\ \cline{1-2} \cline{4-6} 
30 & $L_{3}^{3}$ &  & $y^{[p]}=z$ &  & $p\geq3$ \\ \cline{1-2} \cline{4-6} 
31 & $L_{3}^{4}$ &  & $x^{[p]}=z$&  & $p\geq3$ \\\cline{1-2} \cline{4-6} 
32 & $L_{3}^{5}$ &  & $z^{[p]}=z$ &  & $p\geq3$ \\ \hline
33 & $L_{4}^{1}$ & \multirow{10}{*}{$[w,x]=x$} & $w^{[p]}=w$ &  &  \\ \cline{1-2} \cline{4-6} 
34 & $L_{4}^{2}$ &  & $x^{[p]}=y,~   w^{[p]}=w$ &  &  \\ \cline{1-2} \cline{4-6} 
35 & $L_{4}^{3}$ &  & $y^{[p]}=z,~  w^{[p]}=w$ &  &  \\ \cline{1-2} \cline{4-6} 
36 & $L_{4}^{4}$ &  & $x^{[p]}=z,~ y^{[p]}=z,~  w^{[p]}=w$ &  &  \\ \cline{1-2} \cline{4-6} 
37 & $L_{4}^{5}$ &  & $x^{[p]}=y,~ y^{[p]}=z,~  w^{[p]}=w$ &  &  \\ \cline{1-2} \cline{4-6} 
38 & $L_{4}^{6}$ &  & $y^{[p]}=y,~ w^{[p]}=w$ &  &  \\ \cline{1-2} \cline{4-6}  
39 & $L_{4}^{7}$ &  & $x^{[p]}=z,~ y^{[p]}=y,~  w^{[p]}=w$ &  &  \\ \cline{1-2} \cline{4-6}
40 & $L_{4}^{8}$ &  & $x^{[p]}=y,~ y^{[p]}=y,~  w^{[p]}=w$ &  &  \\ \cline{1-2} \cline{4-6} 
41 & $L_{4}^{9}$ &  & $x^{[p]}=y+z,~ y^{[p]}=y,~  w^{[p]}=w$ &  &  \\ \cline{1-2} \cline{4-6} 
42 & $L_{4}^{10}$ &  & $y^{[p]}=y,~z^{[p]}=z,~  w^{[p]}=w$ &  &  \\ \cline{1-2} \cline{4-6} 
43 & $L_{4}^{11}(\lambda)$ &  & \makecell{$x^{[p]}=y+\lambda z,~   y^{[p]}=y,$ \\ $~ z^{[p]}=z,~ w^{[p]}=w$} & $\lambda\in \mathbb{F}$ &  \\ \hline
44 & $L_{5}^{1}(\xi)$ & \multirow{4}{*}{$[w,x]=x,~[w,y]=\xi y$} & $w^{[p]}=w$ & $\xi \in \Xi_{p}$\tnote{$\dagger$} &  \\ \cline{1-2} \cline{4-6} 
45 & $L_{5}^{2}(\xi)$ &  & $x^{[p]}=z,~   w^{[p]}=w$ & $\xi \in \Xi_{p} $ &  \\ \cline{1-2} \cline{4-6} 
46 & $L_{5}^{3}(\xi)$ &  & $y^{[p]}=z,~ w^{[p]}=w$ & $\xi \in \Xi_{p}\setminus\{\pm1\}$ &  \\ \cline{1-2} \cline{4-6} 
47 & $L_{5}^{4}(\xi)$ &  & $x^{[p]}=z,~ y^{[p]}=z,~  w^{[p]}=w$ & $\xi \in \Xi_{p} \setminus\{1\}$ &  \\ \cline{1-2} \cline{4-6} 
48 & $L_{5}^{5}(\xi)$ & \multirow{4}{*}{$[w,x]=x,~[w,y]=\xi y$} & $z^{[p]}=z,~ w^{[p]}=w$ & $\xi \in \Xi_{p}$ &  \\ \cline{1-2} \cline{4-6} 
49 & $L_{5}^{6}(\xi)$ &  & $x^{[p]}=z,~  z^{[p]}=z,~ w^{[p]}=w$  & $\xi \in \Xi_{p}$ &  \\ \cline{1-2} \cline{4-6} 
50 & $L_{5}^{7}(\xi)$ &  & $y^{[p]}=z,~ z^{[p]}=z,~  w^{[p]}=w$ & $\xi \in \Xi_{p} \setminus\{\pm1\}$ &  \\ \cline{1-2} \cline{4-6} 
51 & $L_{5}^{8}(\xi)$ &  & \makecell{$x^{[p]}=z,~ y^{[p]}=z$,\\ $ z^{[p]}=z,~ w^{[p]}=w$} & $\xi \in \Xi_{p}\setminus\{1\}$ &  \\ \hline
52 & $L_{6}^{1}(\xi,\eta)$ & \makecell{$[w,x]=x,~[w,y]=\xi y$, \\ $ [w,z]=\eta z$} & $w^{[p]}=w$ & $\xi, \eta \in \mathbb{F}_{p}^{\times}$ &  \\ \hline
53 & $N_{1}^{1}$ & $[y,x]=x,~[w,z]=z$ & $y^{[p]}=y,~  w^{[p]}=w$ &  &  \\ \hline
54 & $N_{2}^{1}$ & \multirow{3}{*}{\makecell{$[z,x]=y,~[w,x]=x$,\\ $[w,y]=2y,~[w,z]=z$}} & $w^{[p]}=w$ &  &  \\ \cline{1-2} \cline{4-6} 
55 & $N_{2}^{2}$ &  & $x^{[p]}=y$, $w^{[p]}=w$ & & $p=2$ \\ \cline{1-2} \cline{4-6} 
56 & $N_{2}^{3}$ &  & $y^{[p]}=y$, $w^{[p]}=w$ & & $p=2$ \\ \hline
57 & $N_{3}^{1}(\xi)$ & \multirow{2}{*}{\makecell{$[z,x]=y,~[w,x]=x+\xi z$,\\ $[w,y]=y,~[w,z]=x$}} & $w^{[p]}=w$ & $\xi\in Q_{p}-\frac{1}{4}\tnote{$\ddagger$}$ & $p\geq3$ \\ \cline{1-2} \cline{4-6} 
58 & $N_{3}^{2}(\xi)$ &  & $z^{[p]}=y,~w^{[p]}=w$ & $\xi=0$ & $p=2$ \\ \hline
59 & $N_{4}^{1}$ & \multirow{4}{*}{\makecell{$[z,x]=y,~[w,x]=z$,\\$[w,z]=x$}} & $w^{[p]}=w$ &  & $p\geq3$ \\ \cline{1-2} \cline{4-6} 
60 & $N_{4}^{2}$ &  & $x^{[p]}=y,~   w^{[p]}=w$ &  & $p\geq3$ \\ \cline{1-2} \cline{4-6} 
61 & $N_{4}^{3}$ &  & $y^{[p]}=y,~  w^{[p]}=w$ &  & $p\geq3$ \\ \cline{1-2} \cline{4-6} 
62 & $N_{4}^{4}$ &  & $x^{[p]}=y,~z^{[p]}=y,~  w^{[p]}=w$ &  & $p\geq3$ \\ \hline
63 & $\mathfrak{gl}_{2}^{1}$ & \multirow{5}{*}{\makecell{$[y,x]=-z,~[z,x]=2x,$ \\$[z,y]=-2y$}} & $z^{[p]}=z$ &  & $p\geq3$ \\ \cline{1-2} \cline{4-6} 
64 & $\mathfrak{gl}_{2}^{2}$ &  & $x^{[p]}=w,~z^{[p]}=z$ &  & $p\geq3$ \\ \cline{1-2} \cline{4-6} 
65 & $\mathfrak{gl}_{2}^{3}$ &  & $z^{[p]}=z,~w^{[p]}= w$ &  & $p\geq3$ \\ \cline{1-2} \cline{4-6} 
66 & $\mathfrak{gl}_{2}^{4}$ &  & $x^{[p]}=w,~z^{[p]}=z,~w^{[p]}=w$ &  & $p\geq3$ \\ \cline{1-2} \cline{4-6} 
67 & $\mathfrak{gl}_{2}^{5}(\lambda)$ &  & $z^{[p]}=z+w,~w^{[p]}=\lambda w$ & $\lambda\in\mathbb{F}$ & $p\geq3$ \\ \hline
\end{longtable}
}

\end{ThreePartTable}

\newpage

\section{Codes in Magma}\label{codes}

\begin{lstlisting}[
    style=macaulay,
    caption={\rm Codes for Remark \ref{Groebner}},
    label={autoN},
    numbers=left,
    numberstyle=\tiny\color{gray},
    stepnumber=1,
    numbersep=5pt,
    xleftmargin=15pt
]
> Q:= RationalField();
> P<t,a1,a2,a3,b1,b2,b3,c1,c2,c3> := PolynomialRing(Q, 10);
> I:= ideal<P |
> a1*c3-a3*c1-a1,
> a3*c2-a2*c3-a2,
> a2*c1-a1*c2-2*a3,
> b3*c1-b1*c3-b1,
> b2*c3-b3*c2-b2,
> b1*c2-b2*c1-2*b3,
> 2*a3*b1-2*a1*b3-c1,
> 2*a2*b3-2*a3*b2-c2,
> a1*b2-a2*b1-c3,
> t*(a1*b2*c3+a2*b3*c1+a3*b1*c2-a1*b3*c2-a2*b1*c3-a3*b2*c1)-1>;
> GroebnerBasis(I)
\end{lstlisting}

\begin{lstlisting}[
    style=macaulay,
    caption={\rm Codes for Remark \ref{Groebner-2}},
    label={autoresofN},
    numbers=left,
    numberstyle=\tiny\color{gray},
    stepnumber=1,
    numbersep=5pt,
    xleftmargin=15pt
]
> Q:= RationalField();
> P<t,a1,a2,a3,b1,b2,b3,c1,c2,c3,X,Y> := PolynomialRing(Q, 12);
> J:= ideal<P |
> a1*c3-a3*c1-a1
> a3*c2-a2*c3-a2,
> a2*c1-a1*c2-2*a3,
> b3*c1-b1*c3-b1, 
> b2*c3-b3*c2-b2,
> b1*c2-b2*c1-2*b3,
> 2*a3*b1-2*a1*b3-c1, 
> 2*a2*b3-2*a3*b2-c2,
> a1*b2-a2*b1-c3,
> a3-X,
> b3-X*Y,
> c3-X,
> t*(a1*b2*c3+a2*b3*c1+a3*b1*c2-a1*b3*c2-a2*b1*c3-a3*b2*c1)*X-1>;
> GroebnerBasis(J)
\end{lstlisting}

\Addresses

\end{document}